\documentclass{amsart}

\usepackage{amssymb,mathscinet}
\usepackage{amsthm}
\usepackage{amsmath}
\usepackage[all]{xy} \SelectTips{eu}{}
\usepackage[hidelinks]{hyperref}
\usepackage{mathrsfs}
\usepackage{xfrac,colonequals,rotating}

\newcommand{\numberseries}{\bfseries} 

\newlength{\thmtopspace} 
\newlength{\thmbotspace} 
\newlength{\thmheadspace} 
\newlength{\thmindent} 

\newtheoremstyle{bfupright head,upright body}
{\thmtopspace}{\thmbotspace}
{\upshape}{\thmindent}{\bfseries}{.}{\thmheadspace} {{\numberseries
    \thmnumber{#2\;}}\thmnote{#3}}

\newtheoremstyle{fixed bf head,slanted body}
{\thmtopspace}{\thmbotspace}{\slshape}
{\thmindent}{\bfseries}{.}{\thmheadspace} {{\numberseries
    \thmnumber{#2\;}}\thmname{#1}\thmnote{ (#3)}}

\newtheoremstyle{fixed bf head,upright body}
{\thmtopspace}{\thmbotspace}{\upshape}
{\thmindent}{\bfseries}{.}{\thmheadspace} {{\numberseries
    \thmnumber{#2\;}}\thmname{#1}\thmnote{ (#3)}}

\theoremstyle{bfupright head,upright body} \newtheorem{res}{}[section]

\theoremstyle{fixed bf head,slanted body}
\newtheorem{thm}[res]{Theorem} \newtheorem*{thm*}{Theorem}
\newtheorem{prp}[res]{Proposition} \newtheorem*{prp*}{Proposition}

\newtheorem{cor}[res]{Corollary} \newtheorem*{cor*}{Corollary}
\newtheorem{lem}[res]{Lemma} \newtheorem*{lem*}{Lemma}

\theoremstyle{fixed bf head,upright body}
\newtheorem{dfn}[res]{Definition} \newtheorem*{dfn*}{Definition}
\newtheorem{rmk}[res]{Remark} \newtheorem*{rmk*}{Remark}
\newtheorem{exa}[res]{Example} \newtheorem*{exa*}{Example}
\newtheorem{fct}[res]{Fact} \newtheorem*{fct*}{Fact}

\newlength{\thmlistleft} 
\newlength{\thmlistright} 
\newlength{\thmlistpartopsep} 
\newlength{\thmlisttopsep} 
\newlength{\thmlistparsep} 
\newlength{\thmlistitemsep} 

\newcounter{eqc} \newenvironment{eqc}{\begin{list}{\upshape
      (\textit{\roman{eqc}})}%
    {\usecounter{eqc}%
      \setlength{\leftmargin}{\thmlistleft}%
      \setlength{\labelwidth}{\thmlistleft}%
      \setlength{\rightmargin}{\thmlistright}%
      \setlength{\partopsep}{\thmlistpartopsep}%
      \setlength{\topsep}{\thmlisttopsep}%
      \setlength{\parsep}{\thmlistparsep}%
      \setlength{\itemsep}{\thmlistitemsep}}}%
  {\end{list}}%

\newcommand{\eqclbl}[1]{{\upshape(\textit{#1})}}

\newcounter{prt} \newenvironment{prt}{\begin{list}{\upshape
      (\alph{prt})}%
    {\usecounter{prt}%
      \setlength{\leftmargin}{\thmlistleft}%
      \setlength{\labelwidth}{\thmlistleft}%
      \setlength{\rightmargin}{\thmlistright}%
      \setlength{\partopsep}{\thmlistpartopsep}%
      \setlength{\topsep}{\thmlisttopsep}%
      \setlength{\parsep}{\thmlistparsep}%
      \setlength{\itemsep}{\thmlistitemsep}}}%
  {\end{list}}%

\newcounter{rqm} \newenvironment{rqm}{\begin{list}{\upshape
      (\arabic{rqm})}%
    {\usecounter{rqm}%
      \setlength{\leftmargin}{\thmlistleft}%
      \setlength{\labelwidth}{\thmlistleft}%
      \setlength{\rightmargin}{\thmlistright}%
      \setlength{\partopsep}{\thmlistpartopsep}%
      \setlength{\topsep}{\thmlisttopsep}%
      \setlength{\parsep}{\thmlistparsep}%
      \setlength{\itemsep}{\thmlistitemsep}}}%
  {\end{list}}%

  {\end{list}}%

\newenvironment{prf*}[1][Proof]{%
  \begin{proof}[\bf #1]
    \setcounter{equation}{0}
    } {\end{proof} }

\newcommand{\proofofimp}[3][:]{\mbox{\eqclbl{#2}$\!\implies\!$\eqclbl{#3}#1}}

\newcommand{\thmref}[2][Theorem~]{#1\ref{thm:#2}}
\newcommand{\corref}[2][Corollary~]{#1\ref{cor:#2}}
\newcommand{\prpref}[2][Proposition~]{#1\ref{prp:#2}}
\newcommand{\lemref}[2][Lemma~]{#1\ref{lem:#2}}
\newcommand{\dfnref}[2][Definition~]{#1\ref{dfn:#2}}
\newcommand{\exaref}[2][Example~]{#1\ref{exa:#2}}
\newcommand{\rmkref}[2][Remark~]{#1\ref{rmk:#2}}
\newcommand{\secref}[2][Section~]{#1\ref{sec:#2}}
\newcommand{\fctref}[2][Fact~]{#1\ref{fct:#2}}

\newcommand{\thmcite}[2][?]{\cite[Thm.~#1]{#2}}
\newcommand{\rmkcite}[2][?]{\cite[Rmk.~#1]{#2}}
\newcommand{\exacite}[2][?]{\cite[Exa.~#1]{#2}}
\newcommand{\corcite}[2][?]{\cite[Cor.~#1]{#2}}
\newcommand{\prpcite}[2][?]{\cite[Prop.~#1]{#2}}
\newcommand{\lemcite}[2][?]{\cite[Lem.~#1]{#2}}
\newcommand{\remcite}[2][?]{\cite[Rem.~#1]{#2}}
\newcommand{\seccite}[2][?]{\cite[Sec.~#1]{#2}}
\newcommand{\dfncite}[2][?]{\cite[Def.~#1]{#2}}
\renewcommand{\eqref}[1]{(\ref{eq:#1})} \numberwithin{equation}{res}

\newcommand{\Spec}[1]{\operatorname{Spec}#1}
\newcommand{\qla}{\xla{\qis}}
\newcommand{\qra}{\xra{\qis}}
\newcommand{\Qcoh}[1]{\mathsf{Qcoh}{(#1)}}

\newcommand{\calU}{\mathcal{U}}
\newcommand{\calO}{\mathcal{O}}
\def\urltilda{\kern -.15em\lower .7ex\hbox{\~{}}\kern .04em}
\newcommand{\dif}[2][]{{\partial}^{#2}_{#1}}
\newcommand{\Cone}[1]{\nobreak{\operatorname{Cone}(#1)}}
\newcommand{\cathom}[3]{\operatorname{hom}_{#1}(#2,#3)}

\renewcommand{\leq}{\leqslant}
\renewcommand{\geq}{\geqslant}
\renewcommand{\le}{\leqslant}
\renewcommand{\ge}{\geqslant}
\newcommand{\Ext}[4][A]{\operatorname{Ext}_{#1}^{#2}(#3,#4)}
\newcommand{\ZZ}{\mathbb{Z}}
\newcommand{\Cy}[2][]{\operatorname{Z}_{#1}(#2)}
\newcommand{\Co}[2][]{\operatorname{C}_{#1}(#2)}
\newcommand{\Bo}[2][]{\operatorname{B}_{#1}(#2)}
\newcommand{\lra}{\longrightarrow}
\newcommand{\xla}[2][]{\xleftarrow[#1]{\;#2\;}}
\newcommand{\xra}[2][]{\xrightarrow[#1]{\;#2\;}}
\newcommand{\deq}{\:=\:}

\newcommand{\qqtext}[1]{\qquad\text{#1}\qquad}

\newcommand{\qqand}{\qqtext{and}}
\newcommand{\gra}{\alpha}
\newcommand{\grg}{\gamma}
\newcommand{\grb}{\beta}
\newcommand{\qisdef}[4][\xra{\qis}]{\nobreak{#2\colon #3 #1 #4}}
\newcommand{\mapdef}[4][\rightarrow]{\nobreak{#2\colon #3 #1 #4}}
\newcommand{\dmapdef}[4][\lra]{\nobreak{#2\colon #3\:#1\:#4}}
\newcommand{\is}{\cong}
\newcommand{\qis}{\simeq}
\newcommand{\im}{\operatorname{im}}

\renewcommand{\H}[2][]{\operatorname{H}_{#1}(#2)}
\newcommand{\Thb}[2]{#2_{{\scriptscriptstyle\ge}#1}}
\newcommand{\Tha}[2]{#2_{{\scriptscriptstyle\le}#1}}
\newcommand{\Tsa}[2]{#2_{{\scriptscriptstyle\subseteq}#1}}
\newcommand{\Tsb}[2]{#2_{{\scriptscriptstyle\supseteq}#1}}

\newcommand{\F}{\mathscr{F}}
\newcommand{\C}{\mathscr{C}}
\newcommand{\G}{\mathscr{G}}
\newcommand{\M}{\mathscr{M}}
\newcommand{\semiUV}{semi-$\sfU$-$\sfV$}

\newcommand{\sfA}{\mathsf{A}}

\newcommand{\sfE}{\mathsf{E}}
\newcommand{\sfP}{\mathsf{P}}
\newcommand{\sfS}{\mathsf{S}}
\newcommand{\sfU}{\mathsf{U}}
\newcommand{\sfV}{\mathsf{V}}
\newcommand{\sfUV}{(\sfU,\sfV)}
\newcommand{\capUV}{\sfU\cap\sfV}

\newcommand{\Hom}[3][A]{\operatorname{Hom}_{#1}(#2,#3)}
\newcommand{\fuFPrj}{\operatorname{F_{Prj}}}

\newcommand{\fuFInj}{\operatorname{F_{Inj}}}

\newcommand{\catfinR}[2][\sfU]{\mathsf{Perf_{RGor_{#1}(\sfA)}}(#2)}
\newcommand{\catfinL}[2][\sfV]{\mathsf{Perf_{LGor_{#1}(\sfA)}}(#2)}
\newcommand{\catGL}[2][\sfV]{\mathsf{LGor}_{#1}(#2)}
\newcommand{\catGR}[2][\sfU]{\mathsf{RGor}_{#1}(#2)}
\newcommand{\catG}[2][\sfW]{\mathsf{Gor}_{#1}(#2)}

\newcommand{\catCot}[1]{\mathsf{Cot}(#1)}
\newcommand{\catFlat}[1]{\mathsf{Flat}(#1)}
\newcommand{\catFlatCot}[1]{\mathsf{FlatCot}(#1)}
\newcommand{\catprj}[1]{\mathsf{prj}(#1)}
\newcommand{\catPrj}[1]{\mathsf{Prj}(#1)}
\newcommand{\catInj}[1]{\mathsf{Inj}(#1)}
\newcommand{\catMod}[1]{\mathsf{Mod}(#1)}
\newcommand{\stabcatGR}[2][\sfU]{\mathsf{StRGor}_{#1}(#2)}
\newcommand{\stabcatGL}[2][\sfV]{\mathsf{StLGor}_{#1}(#2)}
\newcommand{\stabcatG}[2][\sfW]{\mathsf{StGor}_{#1}(#2)}
\newcommand{\sfFlat}{\mathsf{Flat}}
\newcommand{\sfprj}{\mathsf{prj}}
\newcommand{\sfPrj}{\mathsf{Prj}}
\newcommand{\sfFlatCot}{\mathsf{FlatCot}}
\newcommand{\catDFtac}[1]{\mathsf{D}_{\operatorname{F-tac}} (#1)}

\newcommand{\dle}{\:\le\:}

\newcommand{\dqis}{\:\qis\:}

\newcommand{\catK}[1]{\mathsf{K}(#1)}
\newcommand{\catD}[1]{\mathsf{D}(#1)}
\newcommand{\catDsg}[1]{\mathsf{D}_{\operatorname{sg}}(#1)}
\newcommand{\catDsgInj}[1]{\mathsf{D}^{\operatorname{sg}}_{\operatorname{Inj}}(#1)}
\newcommand{\catDsgPrj}[1]{\mathsf{D}^{\operatorname{sg}}_{\operatorname{Prj}}(#1)}
\newcommand{\catDb}[1]{\mathsf{D}_\mathrm{b}(#1)}

\newcommand{\catDdefInj}[1]{\mathsf{D}^{\operatorname{def}}_{\operatorname{Inj}}(#1)}
\newcommand{\catDdefPrj}[1]{\mathsf{D}^{\operatorname{def}}_{\operatorname{Prj}}(#1)}
\newcommand{\catKb}[1]{\mathsf{K}_\mathrm{b}(#1)}
\newcommand{\catKFtac}[1]{\mathsf{K}_{\operatorname{F-tac}}(#1)}
\newcommand{\catKRtac}[2][\sfU]{\mathsf{K}^{\operatorname{R}}_{\operatorname{#1-tac}}(#2)}
\newcommand{\catKLtac}[2][\sfV]{\mathsf{K}^{\operatorname{L}}_{\operatorname{#1-tac}}(#2)}
\newcommand{\catKac}[1]{\mathsf{K}_{\operatorname{ac}}(#1)}
\newcommand{\catKeac}[2][\sfE]{\mathsf{K}_{\operatorname{#1-ac}}(#2)}
\newcommand{\catKeeac}[2][\sfE]{\mathsf{K}_{\operatorname{(#1-ac)}}(#2)}
\newcommand{\catKpac}[1]{\mathsf{K}_{\operatorname{pac}}(#1)}
 \newcommand{\catKtac}[1]{\mathsf{K}_{\operatorname{tac}}(#1)}
\newcommand{\Ftac}{{\bf F}-totally acyclic}
\newcommand{\catC}[1]{\mathsf{C}(#1)}
\newcommand{\bcatC}[2]{\mathsf{C}_{\textnormal{#1}}(#2)}
\newcommand{\bcatK}[2]{\mathsf{K}_{\textnormal{#1}}(#2)}
\newcommand{\bcatD}[2]{\mathsf{D}_{\textnormal{#1}}(#2)}

\newcommand{\bcatKeac}[3][\sfE]{\mathsf{K}_{\operatorname{#1-ac},\textnormal{#2}}(#3)}
\newcommand{\bcatKeeac}[3][\sfE]{\mathsf{K}_{\operatorname{(#1-ac)},\textnormal{#2}}(#3)}

\newcommand{\Shift}[2][]{\Sigma^{#1}{#2}}
\newcommand{\findim}[1]{\operatorname{findim}(#1)}
\newcommand{\JSt}{\v{S}\soft{t}ov\'{\i}\v{c}ek}

\usepackage{xcolor}

\def\soft#1{\leavevmode\setbox0=\hbox{h}\dimen7=\ht0\advance \dimen7
  by-1ex\relax\if t#1\relax\rlap{\raise.6\dimen7
  \hbox{\kern.3ex\char'47}}#1\relax\else\if T#1\relax
  \rlap{\raise.5\dimen7\hbox{\kern1.3ex\char'47}}#1\relax \else\if
  d#1\relax\rlap{\raise.5\dimen7\hbox{\kern.9ex \char'47}}#1\relax\else\if
  D#1\relax\rlap{\raise.5\dimen7 \hbox{\kern1.4ex\char'47}}#1\relax\else\if
  l#1\relax \rlap{\raise.5\dimen7\hbox{\kern.4ex\char'47}}#1\relax \else\if
  L#1\relax\rlap{\raise.5\dimen7\hbox{\kern.7ex
  \char'47}}#1\relax\else\message{accent \string\soft \space #1 not
  defined!}#1\relax\fi\fi\fi\fi\fi\fi}

\title[The singularity category of an exact category]{The singularity
  category of an exact category applied to characterize Gorenstein schemes}

\author[L.W.\ Christensen]{Lars Winther Christensen} \address{L.W.C. \
  Texas Tech University, Lubbock, TX 79409, U.S.A.}
\email{lars.w.christensen@ttu.edu}
\urladdr{http://www.math.ttu.edu/\urltilda lchriste}

\author[N.\ Ding]{Nanqing Ding} \address{N.D. \ Nanjing University,
  Nanjing 210093, China} \email{nqding@nju.edu.cn} \urladdr{}

\author[S.\ Estrada]{Sergio Estrada} \address{S.E. \ Universidad de
  Murcia, Murcia 30100, Spain} \email{sestrada@um.es}
\urladdr{https://webs.um.es/sestrada/}

\author[J.\ Hu]{Jiangsheng Hu} \address{J.H. \ Jiangsu University of
  Technology, Changzhou 213001, China} \email{jiangshenghu@hotmail.com} \urladdr{}

\author[H.\ Li]{Huanhuan Li} \address{H.L. \ Xidian University, Xi'an 710071,
    China} \email{lihh@xidian.edu.cn}

\author[P.\ Thompson]{Peder Thompson} \address{P.T. \ Norwegian
  University of Science and Technology, 7491 Trondheim, Norway}
\email{peder.thompson@ntnu.no} \urladdr{http://pthompson.nupurple.net/}
\curraddr{Niagara University, Niagara University, NY 14109, U.S.A.}

\thanks{L.W.C.\ was partly supported by Simons Foundation grant
  428308. N.D. was supported by NSFC grant 11771202. S.E. was partly
  supported by grant PID2020-113206GB-I00 funded by
  MCIN/AEI/10.13039/501100011033.  J.H. was partly supported by NSFC
  grant 11771212. H.L. was partly supported by NSFC grant
  11626179. Part of the work was done when J.H. visited Universidad de
  Murcia in February 2020; the hospitality of this institution is
  acknowledged with gratitude.}

\date{5 May 2022}



\begin{document}

\begin{abstract}
  We construct a non-affine analogue of the singularity category of a
  Gorenstein local ring. With this Buchweitz's classic equivalence of
  three triangulated categories over a Gorenstein local ring has been
  generalized to schemes, a project started by Murfet and Salarian
  more than ten years ago. Our construction originates in a framework
  we develop for singularity categories of exact categories. As an
  application of this framework in the non-commutative setting, we
  characterize rings of finite finitistic dimension.
\end{abstract}


\maketitle
\thispagestyle{empty}

\section*{Introduction}

\noindent
The singularity category of a commutative noetherian ring $A$ is the
Verdier quotient of the finite bounded derived category of $A$ by the
subcategory of perfect $A$-complexes. A motivation for this paper
comes from work of Buchweitz \cite{ROB86}: In contemporary
terminology, he proved that the singularity category of a Gorenstein
local ring $A$ is equivalent to the stable category of finitely
generated Gorenstein projective $A$-modules and to the homotopy
category of totally acyclic complexes of finitely generated projective
$A$-modules; in symbols
\begin{equation}\tag{$\diamond$}
  \catDsg{A} \dqis \stabcatG[\sfprj]{A} \dqis \catKtac{\catprj{A}}  \:.
\end{equation}
A perfect non-affine analogue may be beyond reach: Work of Orlov
\cite{DOO04} shows that even under the assumption that a scheme has
enough locally free sheaves---a prerequisite for talking about
Gorenstein projective coherent sheaves---these sheaves do not form a
Frobenius category; i.e.\ they do not provide an analogue of
$\stabcatG[\sfprj]{A}$.

Murfet and Salarian \cite{DMfSSl11} initiated the quest for non-affine
analogues of the categories in $(\diamond)$. Noticing that the
structure sheaf of a scheme is flat, but not necessarily projective,
their idea is to replace projective objects by flat objects. In this
process, though, one must give up finite generation. Fortunately,
there are ``big'' versions of all three categories in
$(\diamond)$---obtained by dropping the assumptions of finite
generation---and they are also equivalent. Indeed, the big singularity
category of a Gorenstein ring $A$ is by work of Beligiannis
\thmcite[6.9]{ABl00} equivalent to the stable category
$\stabcatG[\sfPrj]{A}$ of Gorenstein projective $A$-modules, and the
latter category is for every ring $A$ equivalent to the homotopy
category $\catKtac{\catPrj{A}}$ of totally acyclic complexes of
projective $A$-modules, see for example \corcite[3.9]{CET-20}. Thus,
these are equivalences one can aim to replicate in the non-affine
setting.

As an analogue of the category $\catKtac{\catPrj{A}}$ for a
semi-separated noetherian scheme $X$, Murfet and Salarian identified
the category
\begin{equation*}
  \catDFtac{\catFlat{X}}
  \,\colonequals\,
  \frac{\catKFtac{\catFlat{X}}}{\catKpac{\catFlat{X}}} \:,
\end{equation*}
i.e.\ the Verdier quotient of the homotopy category of \Ftac\
complexes of flat quasi-coherent sheaves on $X$ by its subcategory of
pure-acyclic complexes.  Indeed, for a commutative noetherian ring $A$
of finite Krull dimension and the scheme $X=\Spec(A)$, the categories
$\catKtac{\catPrj{A}}$ and $\catDFtac{\catFlat{X}}$ are equivalent by
\lemcite[4.22]{DMfSSl11}. An analogue for schemes of the stable
category was identified in \cite{CET}, based on a change of focus from
flat objects to flat-cotorsion objects: For a semi-separated
noetherian scheme $X$, the category $\catDFtac{\catFlat{X}}$ is
equivalent to the homotopy category $\catKtac{\catFlatCot{X}}$ of
totally acyclic complexes of flat-cotorsion sheaves on $X$, and that
category is equivalent to the stable category
$\stabcatG[\sfFlatCot]{X}$ of Gorenstein flat-cotorsion sheaves on
$X$. The primary goal of this paper is to identify an analogue of the
singularity category in the non-affine setting and to show that it is
equivalent to the categories $\catKtac{\catFlatCot{X}}$ and
$\stabcatG[\sfFlatCot]{X}$ for a Gorenstein scheme $X$ of finite Krull
dimension.

We achieve this goal as an application of a more general theory that
we develop for singularity categories associated to a cotorsion pair.
Before giving an outline of this theory, we remark that in the case of
flat-cotorsion sheaves on $X$, the singularity category takes the form
\begin{equation*}
  \frac{\catDb{\catCot{X}}}
  {\langle\text{bounded complexes of flat-cotorsion sheaves}\rangle} \:.
\end{equation*}
We notice in \exaref{Regsingcat} that this singularity category
detects when a noetherian semi-separated scheme of finite Krull
dimension is regular. Further, for a Gorenstein scheme of finite Krull
dimension the category is compactly generated and closely related to
Orlov's singularity category, see \rmkref{gorX}.

Here is the outline: Let $\sfA$ be an abelian category. To an additive
full subcategory $\sfE$ of $\sfA$ that is closed under extensions and
direct summands, we associate in \secref{exact} two singularity
categories relative to the projective and injective objects in
$\sfE$. If $\sfE$ is part of a complete hereditary cotorsion pair,
then there are natural functors from categories of Gorenstein objects
into these singularity categories. This allows us in
\secref{defect}---inspired by works of Bergh, Jorgensen, and Oppermann
\cite{BJO-15} and Iyengar and Krause \cite{SInHKr06}---to define
associated defect categories. Given a complete hereditary cotorsion
pair $\sfUV$ in $\sfA$, we build on the theory from \cite{CET-20} to
develop Gorenstein dimensions for objects in $\sfU$ and
$\sfV$. Assuming that $\sfA$ is Grothendieck---as is the case for the
categories of modules over a ring and quasi-coherent sheaves on a
semi-separated noetherian scheme---these dimensions extend through
work of Gillespie to invariants on the derived category of $\sfA$;
that's the topic of \secref{gdim}.  In \secref{appl1} we arrive at
equivalences of categories akin to those in $(\diamond)$---but now for
a Gorenstein scheme. In \secref{appl2} we characterize rings of finite
finitistic dimension in terms of the singularity and defect categories
developed in the first sections.
\begin{equation*}
  * \ \ * \ \ *
\end{equation*}
\noindent
Throughout the paper, $\sfA$ denotes an abelian category; hom-sets of
objects in $\sfA$ are denoted $\operatorname{hom}_\sfA$, and that
notation is also used for the induced functor from
$\sfA^\mathrm{op} \times \sfA$ to abelian groups. By a
\emph{subcategory} of $\sfA$ we mean a full subcategory that is closed
under isomorphisms.  Let $\sfS$ be a subcategory of $\sfA$. The
category of chain complexes of objects from $\sfS$, or
$\sfS$-complexes, is denoted $\catC{\sfS}$ and $\catK{\sfS}$ is the
associated homotopy category; the latter is triangulated if $\sfS$ is
additive. A complex $X$ in $\catC{\sfS}$ is said to be \emph{bounded
  below} if $X_n = 0$ holds for $n \ll 0$, \emph{bounded above} if
$X_n = 0$ holds for $n \gg 0$, and \emph{bounded} if $X_n = 0$ holds
for $|n| \gg 0$.  The full subcategories of bounded below, bounded
above, and bounded complexes are denoted $\bcatC{+}{\sfS}$,
$\bcatC{$-$}{\sfS}$, and $\bcatC{b}{\sfS}$, respectively. The
essential images of these subcategories in the homotopy category
$\catK{\sfS}$ are triangulated subcategories denoted
$\bcatK{+}{\sfS}$, $\bcatK{$-$}{\sfS}$, and $\bcatK{b}{\sfS}$,
respectively.

Given a complex $X \in \catC{\sfA}$ with differential $\dif{X}$ we
write $\Bo[n]{X}$ and $\Cy[n]{X}$ for the subobjects of boundaries and
cycles in $X_n$. The cokernel of $\dif[n+1]{X}$ is the quotient object
$\Co[n]{X} \colonequals X_n/\Bo[n]{X}$, and the homology of $X$ in
degree $n$ is the subquotient object
$\H[n]{X} \colonequals \Cy[n]{X}/\Bo[n]{X}$. Given a second
$\sfA$-complex $Y$, we write $\Hom[\sfA]{X}{Y}$ for the total
hom-complex in $\catC{\sfA}$.

\section{Singularity categories of exact categories}
\label{sec:exact}

\noindent
Let $\sfE \subseteq \sfA$ be an additive subcategory closed under
extensions and direct summands; it is an idempotent complete exact
category with the exact structure induced from $\sfA$, see
B\"uhler~\rmkcite[6.2 and Lem.~10.20]{TBh10}. The setup developed in
this first section could be done in the generality of idempotent
complete exact categories. We don't need that for the applications we
pursue in this paper, but we structure the arguments in such a way as
to make the generalizations straightforward for the initiated reader.

Given a subcategory $\sfS \subseteq \sfE$, a complex $X$ in
$\catK{\sfS}$ is called \emph{$\sfE$-acyclic} if
\begin{equation}
  \label{eq:ZXZ}
  0 \lra \Cy[n]{X} \lra X_n \lra \Cy[n-1]{X} \lra 0
\end{equation}
is an exact sequence in $\sfE$ for every $n\in\ZZ$. The full
subcategory of such complexes is denoted $\catKeac{\sfS}$; it is
isomorphism closed, see Keller \cite[4.1]{BKl90}, and is indeed
triangulated, see Neeman \lemcite[1.1]{ANm90}. The symbol
$\catKeac{\sfS}$ is used with subscripts $+$, $-$, and $\mathrm{b}$ to
indicate boundedness:
$\bcatKeac{b}{\sfS} \colonequals \catKeac{\sfS} \cap \bcatK{b}{\sfS}$
etc. An $\sfA$-acyclic complex is simply referred to as
\emph{acyclic}, and the corresponding category is denoted
$\catKac{\sfS}$. A morphism in $\catK{\sfS}$ is called an
\emph{$\sfE$-quasi-isomorphism} if its mapping cone is $\sfE$-acyclic
and simply a \emph{quasi-isomorphism} if its mapping cone is acyclic.

The \emph{derived category} of $\sfE$ is defined to be the Verdier
quotient
\begin{equation*}
  \catD{\sfE} \colonequals
  \catK{\sfE}/\catKeac{\sfE} \:;
\end{equation*}
see \cite{ANm90} and Keller \cite{BKl96}.  In the case where
$\sfE = \sfA$ is the module category of a ring, this construction
yields the usual derived category.  The bounded versions
$\bcatD{+}{\sfE}$, $\bcatD{$-$}{\sfE}$, and $\bcatD{b}{\sfE}$ are
defined analogously:
$\bcatD{+}{\sfE} \colonequals \bcatK{+}{\sfE}/\bcatKeac{+}{\sfE}$ etc.

A complex $X$ in $\catK{\sfS}$ is called \emph{eventually
  $\sfE$-acyclic} if the sequences \eqref{ZXZ} are exact in $\sfE$ for
all $|n| \gg 0$. The full subcategory of such complexes is denoted
$\catKeeac{\sfS}$; this notation is also used with subscripts to
indicate boundedness.

Recall that for an $\sfA$-complex $X$ and an integer $n$, the
\emph{hard truncation below of $X$ at $n$} is the complex
$\Thb{n}{X} \colonequals \cdots \to X_{n+1} \to X_n\to 0$, and the
\emph{soft truncation below of $X$ at $n$} is the complex
$\Tsb{n}{X} \colonequals \cdots \to X_{n+1} \to \Cy[n]{X} \to 0$.
Hard and soft truncations above are defined similarly:
$\Tha{n}{X}\colonequals 0\to X_n\to X_{n-1}\to \cdots$ and
$\Tsa{n}{X}\colonequals 0 \to \Co[n]{X}\to X_{n-1} \to \cdots$.

\begin{prp}
  \label{prp:K-ev-ac}
  Let $\sfE$ be an additive subcategory of $\sfA$ that is closed under
  extensions and direct summands; let $\sfS$ be an additive
  subcategory of $\sfE$.
  \begin{prt}
  \item The category $\bcatKeeac{+}{\sfS}$ is a triangulated
    subcategory of $\bcatK{+}{\sfS}$.
  \item The category $\bcatKeeac{$-$}{\sfS}$ is a triangulated
    subcategory of $\bcatK{$-$}{\sfS}$.
  \end{prt}
\end{prp}

\begin{prf*}
  We prove (a); the proof of (b) is similar.  That
  $\bcatKeeac{+}{\sfS}$ is additive and closed under shifts is
  evident; it remains to show that it is closed under isomorphisms and
  cones.

  Let $\mapdef{\gra}{X}{Y}$ be an isomorphism in $\bcatK{+}{\sfS}$
  with $X \in \bcatKeeac{+}{\sfS}$. The complex $\Cone{\gra}$ is
  contractible and hence $\sfE$-acyclic as $\sfE$ is closed under
  direct summands.  For every $n \in \ZZ$ the short exact sequence
  \begin{equation}
    \label{eq:0}
    0 \lra Y_n \lra \Cone{\gra}_n \lra (\Shift{X})_n \lra 0
  \end{equation}
  yields a commutative diagram
  \begin{equation}
    \label{eq:1}
    \begin{gathered}
      \xymatrix{
        \Cone{\gra}_n \ar[r] \ar[d]& X_{n-1}\ar[d] \ar[r] & 0\\
        \Cy[n-1]{\Cone{\gra}} \ar[r] & \Cy[n-2]{X} \:. }
    \end{gathered}
  \end{equation}
  For all $n \gg 0$ the vertical morphisms in \eqref{1} are
  epimorphisms, hence so is the lower horizontal morphism. Thus the
  vertical morphisms in the commutative diagram
  \begin{equation*}
    \xymatrix{
      0 \ar[r] & \Cy[n]{\Cone{\gra}} \ar[r] \ar[d]
      & \Cone{\gra}_{n} \ar[r] \ar[d]
      & \Cy[n-1]{\Cone{\gra}} \ar[r] \ar[d] & 0
      \\
      0 \ar[r] & \Cy[n-1]{X} \ar[r] & X_{n-1} \ar[r] & \Cy[n-2]{X} \ar[r] & 0
    }
  \end{equation*}
  are epimorphisms. As the cycle functor is left exact, \eqref{0} and
  the Snake Lemma yield the exact sequence
  \begin{equation*}
    0 \lra \Cy[n]{Y} \lra Y_n \lra \Cy[n-1]{Y} \lra 0 \:.
  \end{equation*}
  In \eqref{1} the upper horizontal and right-hand vertical morphisms
  are epimorphisms with kernels in $\sfE$, namely $Y_n$ and
  $\Cy[n-1]{X}$. As $\sfE$ is an exact category closed under direct
  summands, \prpcite[7.6]{TBh10} shows that the lower horizontal
  morphism in \eqref{1} is an epimorphism with kernel in $\sfE$. By
  left exactness of the cycle functor, it now follows from \eqref{0}
  that $\Cy[n-1]{Y}$ is in $\sfE$; thus $Y$ belongs to
  $\bcatKeeac{+}{\sfS}$.

  Finally, if $\mapdef{\gra}{X}{Y}$ is a morphism in
  $\bcatKeeac{+}{\sfS}$, then for $n \gg 0$ the complexes $\Tsb{n}{X}$
  and $\Tsb{n}{Y}$ are $\sfE$-acyclic and $\gra$ induces a map $\gra'$
  between them. In high degrees, the cones of $\gra$ and $\gra'$
  agree, and the cone of $\gra'$ is $\sfE$-acyclic; see
  \lemcite[10.3]{TBh10}.
\end{prf*}

We denote by $\catPrj{\sfE}$ and $\catInj{\sfE}$ the subcategories of
projective and injective objects in $\sfE$.

\begin{prp}
  \label{prp:semi-exact}
  Let $\sfE$ be an additive subcategory of $\sfA$ that is closed under
  extensions and direct summands.
  \begin{prt}
  \item If\, $\sfE$ has enough projectives, then there are
    triangulated equivalences of triangulated categories
    \begin{equation*}
      \bcatK{+}{\catPrj{\sfE}} \dqis \bcatD{+}{\sfE}  \qqand
      \bcatKeeac{+}{\catPrj{\sfE}} \dqis \bcatD{b}{\sfE}  \:.
    \end{equation*}

  \item If\, $\sfE$ has enough injectives, then there are triangulated
    equivalences of triangulated categories
    \begin{equation*}
      \bcatK{$-$}{\catInj{\sfE}}  \dqis \bcatD{$-$}{\sfE} \qqand
      \bcatKeeac{$-$}{\catInj{\sfE}}  \dqis \bcatD{b}{\sfE} \:.
    \end{equation*}
  \end{prt}
\end{prp}

\begin{prf*}
  We prove part (a); the proof of part (b) is dual. For every complex
  $X$ in $\bcatK{+}{\sfE}$ there is a distinguished triangle in
  $\bcatK{+}{\sfE}$,
  \begin{equation*}
    p(X) \lra X \lra a(X) \lra \Shift{p(X)} \:,
  \end{equation*}
  with $p(X)$ in $\bcatK{+}{\catPrj{\sfE}}$ and $a(X)$ an
  $\sfE$-acyclic complex; moreover the inclusion of $\sfE$-acyclic
  complexes into $\bcatK{+}{\sfE}$ has a left adjoint.  For a proof
  (of the dual result) see \lemcite[4.1]{BKl90}. It is now standard,
  see Krause \prpcite[4.9.1 and the proof of Prop.\ 4.13.1]{HKr10},
  that the functor
  $G \colon \bcatK{+}{\catPrj{\sfE}} \to \bcatK{+}{\sfE} \to
  \bcatD{+}{\sfE}$ yields the equivalence
  $ \bcatK{+}{\catPrj{\sfE}} \qis \bcatD{+}{\sfE}$; see also
  \exacite[12.2]{BKl96}.

  By \prpref{K-ev-ac} the triangulated subcategory
  $\bcatKeeac{+}{\catPrj{\sfE}}$ is equivalent to its essential image
  under $G$ in $\bcatD{+}{\sfE}$. For $X$ in $\bcatK{b}{\sfE}$ the
  complex $p(X)$ belongs to $\bcatKeeac{+}{\catPrj{\sfE}}$, and a
  complex in $\bcatKeeac{+}{\catPrj{\sfE}}$ is $\sfE$-quasi-isomorphic
  in $\bcatK{+}{\sfE}$ to a complex in $\bcatK{b}{\sfE}$. As
  $\bcatD{b}{\sfE}$ is equivalent to the subcategory of $\catD{\sfE}$
  generated by $\bcatK{b}{\sfE}$, see \lemcite[11.7]{BKl96}, the
  restriction of $G$ to $\bcatKeeac{+}{\catPrj{\sfE}}$ induces the
  second equivalence.
\end{prf*}

As $\bcatK{b}{\catPrj{\sfE}}$ evidently is a thick subcategory of
$\bcatK{+}{\catPrj{\sfE}}$, it follows that for a category $\sfE$ as
in \prpref{semi-exact} with enough projectives, the subcategory
$\bcatK{b}{\catPrj{\sfE}}$ is equivalent to a thick subcategory of
$\bcatD{b}{\sfE}$. Similarly, if $\sfE$ has enough injectives, then
$\bcatK{b}{\catInj{\sfE}}$ is equivalent to a thick subcategory of
$\bcatD{b}{\sfE}$.

\begin{dfn}
  \label{dfn:Dsg}
  Let $\sfE$ be an additive subcategory of $\sfA$ that is closed under
  extensions and direct summands.
  \begin{prt}
  \item Assume that $\sfE$ has enough projectives. The Verdier
    quotient
    \begin{equation*}
      \catDsgPrj{\sfE} \colonequals \catDb{\sfE}/\bcatK{b}{\catPrj{\sfE}}
    \end{equation*}
    is called the \emph{projective singularity category of $\sfE$} or
    $\catPrj{\sfE}$-\emph{singularity category.}

  \item Assume that $\sfE$ has enough injectives. The Verdier quotient
    \begin{equation*}
      \catDsgInj{\sfE} \colonequals \catDb{\sfE}/\bcatK{b}{\catInj{\sfE}}
    \end{equation*}
    is called the \emph{injective singularity category of $\sfE$} or
    $\catInj{\sfE}$-\emph{singularity category.}
  \end{prt}
\end{dfn}

With an eye towards module categories, the next result justifies the
term ``singularity category.'' Resolutions of objects in exact
categories are defined in the usual way, see \seccite[12]{TBh10}.

\begin{thm}
  \label{thm:findim1}
  Let $\sfE$ be an additive subcategory of $\sfA$ that is closed under
  extensions and direct summands.
  \begin{prt}
  \item Assume that $\sfE$ has enough projectives. The category
    $\catDsgPrj{\sfE}$ vanishes if and only if every object in $\sfE$
    has a bounded resolution by objects from $\catPrj{\sfE}$.
  \item Assume that $\sfE$ has enough injectives. The category
    $\catDsgInj{\sfE}$ vanishes if and only if every object in $\sfE$
    has a bounded coresolution by objects from $\catInj{\sfE}$.
  \end{prt}
\end{thm}

\begin{prf*}
  We prove (a); the proof of (b) is dual.  The category
  $\catDsgPrj{\sfE}$ is per \prpref{semi-exact} equivalent to
  $\bcatKeeac{+}{\catPrj{\sfE}}/\bcatK{b}{\catPrj{\sfE}}$. The ``only
  if'' statement is clear. For the ``if'' statement let $P$ be a
  complex in $\bcatKeeac{+}{\catPrj{\sfE}}$. For $n \gg 0$ the objects
  $\Cy[n]{P}$ belong to $\sfE$, and the sequences
  $0 \to \Cy[n+1]{P} \to P_{n+1} \to \Cy[n]{P} \to 0$ are exact. From
  the Comparison Theorem for projective resolutions
  \thmcite[12.4]{TBh10} it follows that $P$ is isomorphic to a complex
  in $\bcatK{b}{\catPrj{\sfE}}$.
\end{prf*}

Given a ring $A$ we write $\catMod{A}$ for the abelian category of
$A$-modules.

\begin{exa}
  Let $A$ be a ring. The category $\catDsgPrj{\catMod{A}}$ vanishes if
  and only if every $A$-module has finite projective dimension; that
  is, if and only if $A$ has finite global dimension. In particular,
  for a commutative noetherian local ring $A$ the category
  $\catDsgPrj{\catMod{A}}$ provides a measure of how singular---i.e.\
  how far from being regular---$A$ is. It thus serves the same purpose
  as the classic singularity category $\catDsg{A}$ recalled in the
  introduction.
\end{exa}

\begin{rmk}
  A variety of relative singularity categories are considered in the
  literature. They are in certain aspects more general and/or more
  restrictive than our notions. Closest to ours is, perhaps, the
  category considered by Chen~\cite{XWC11}: In the case $\sfE=\sfA$
  his category $D_{\catPrj{\sfA}}(\sfA)$ agrees with
  $\catDsgPrj{\sfA}$.
\end{rmk}

Let $(\sfU,\sfV)$ be a cotorsion pair in $\sfA$; the next result shows
that the setup developed above applies to $\sfU$ and $\sfV$, and in
the balance of the paper, we work in that~setting.  A cotorsion pair
$(\sfU,\sfV)$ is called \emph{complete} if for each object $M\in \sfA$
there are exact sequences $0\to V\to U \to M \to 0$ and
$0 \to M \to V'\to U'\to 0$ with $U,U'\in \sfU$ and $V,V'\in
\sfV$. Such sequences are known as special $\sfU$-precovers and
special $\sfV$-envelopes, respectively. A cotorsion pair is called
\emph{hereditary} if $\Ext[\sfA]{i}{U}{V}=0$ holds for all objects
$U \in \sfU$ and $V \in \sfV$ and all $i >0$.

\begin{prp}
  \label{prp:UV-exact}
  Let $(\sfU,\sfV)$ be a complete cotorsion pair in $\sfA$.
  \begin{prt}
  \item The category $\sfV$ has enough projectives, and one has
    $\catPrj{\sfV} = \capUV$.

  \item The category $\sfU$ has enough injectives, and one has
    $\catInj{\sfU} = \capUV$.
  \end{prt}
\end{prp}

\begin{prf*}
  We prove (a); the proof of (b) is dual. Let $M \in \sfV$; since
  $(\sfU, \sfV)$ is a complete cotorsion pair in $\sfA$, there is an
  exact sequence
  \begin{equation*}
    0 \lra V \lra U \lra M \lra 0
  \end{equation*}
  with $V \in \sfV$ and $U \in \sfU$. As $\sfV$ is closed under
  extensions, the object $U$ is in $\capUV$. It is evident that
  objects in $\capUV$ are projective in $\sfV$, and given any
  projective object $M$ in $\sfV$ the exact sequence above shows that
  it is a quotient object, and hence a direct summand, of an object in
  $\capUV$.  Thus one has $\catPrj{\sfV} = \capUV$.
\end{prf*}

\begin{rmk}
  \label{rmk:findim2}
  Let $(\sfU,\sfV)$ be a complete cotorsion pair.  It follows in view
  of \thmref{findim1} that the projective singularity category of
  $\sfV$ vanishes if and only if every object in $\sfA$ has a bounded
  resolution by objects from $\sfU$. Indeed, every object in $\sfA$
  has a resolution by objects from $\sfU$ in which the first syzygy
  belongs to $\sfV$, and $\catDsgPrj{\sfV}$ vanishes if and only if
  this object has a bounded resolution by objects from
  $\capUV$. Similarly, vanishing of the category $\catDsgInj{\sfU}$
  means that every object in $\sfA$ has a bounded coresolution by
  objects from $\sfV$.

  As $\sfU$ contains $\catPrj{\sfA}$ the singularity category
  $\catDsgPrj{\sfU}$ is defined if $\sfA$ has enough projectives, but
  we are not aware of any interpretation of it in this generality. For
  the absolute cotorsion pair $(\catPrj{\sfA},\sfA)$ in a category
  with enough projectives, the singularity category
  $\catDsgPrj{\catPrj{\sfA}}$, of course, vanishes. The category
  $\catDsgInj{\sfV}$ is defined if $\sfA$ has enough injectives and in
  general it appears equally intractable.
\end{rmk}

\begin{exa}
  \label{exa:FCsingcat}
  Let $A$ be a ring and consider the cotorsion pair
  $(\catFlat{A},\catCot{A})$, which is complete.

  The singularity category $\catDsgPrj{\catCot{A}}$ vanishes if and
  only if every $A$-module has finite flat dimension; that is, if and
  only if $A$ has finite weak global dimension.

  The singularity category $\catDsgInj{\catFlat{A}}$ vanishes if and
  only if every $A$-module has finite cotorsion dimension; by a result
  of Mao and Ding \thmcite[19.2.5]{LMaNDn06} this is equivalent to $A$
  being $n$-perfect---every flat module has projective dimension at
  most $n$---for some $n$.

  If $A$ has finite global dimension, then both
  $\catDsgPrj{\catCot{A}}$ and $\catDsgInj{\catFlat{A}}$ vanish, and
  the converse holds by \thmcite[19.2.14]{LMaNDn06}.
\end{exa}

\begin{exa}
  \label{exa:Regsingcat}
  Let $X$ be a semi-separated noetherian scheme of finite Krull
  dimension. The flat sheaves and cotorsion sheaves on $X$ form a
  complete cotorsion pair $(\catFlat{X},\catCot{X})$ in the category
  of quasi-coherent sheaves on $X$; see \rmkcite[2.4]{CET}. Now the
  singularity category $\catDsgPrj{\catCot{X}}$ vanishes if and only
  if every sheaf has finite flat dimension. This is equivalent to
  ${\mathcal O}_{X,x}$ being a regular local ring for every $x\in X$;
  that is $\catDsgPrj{\catCot{X}}$ vanishes if and only if $X$ is a
  regular scheme.
\end{exa}

\section{Gorenstein defect categories}
\label{sec:defect}

\noindent
In this section, $(\sfU,\sfV)$ is a complete cotorsion pair in $\sfA$.
The additive categories $\sfU$ and $\sfV$ are closed under extensions
and direct summands, and by \prpref{UV-exact} the category $\sfV$ has
enough projectives, and $\sfU$ has enough injectives. Thus the
singularity categories $\catDsgPrj{\sfV}$ and $\catDsgInj{\sfU}$ are
defined in this context. In the case of a complete hereditary
cotorsion pair we give a more concrete interpretation of these
singularity categories in terms of the notions of $\sfU$- and
$\sfV$-Gorenstein objects.

We recall from \cite{CET-20} that an acyclic complex $T$ of objects
from $\sfU$ is called \emph{right $\sfU$-totally acyclic} if and only
if the cycle objects $\Cy[n]{T}$ belong to $\sfV$ and the complex
$\cathom{\sfA}{T}{W}$ is acyclic for every $W \in \capUV$. As $\sfV$
is extension closed, the objects in a $\sfU$-totally acyclic complex
in fact belong to $\capUV$; the symbol $\catKRtac{\capUV}$ denotes the
homotopy category of such complexes.  An object is called \emph{right
  $\sfU$-Gorenstein} if it equals $\Cy[0]{T}$ for some right
$\sfU$-totally acyclic complex $T$.  The subcategory of right
$\sfU$-Gorenstein objects in $\sfA$ is denoted $\catGR{\sfA}$; it is
by \thmcite[2.11]{CET-20} a Frobenius category, and the associated
stable category is denoted $\stabcatGR{\sfA}$. Dually one defines
\emph{left $\sfV$-totally acyclic} complexes, \emph{left
  $\sfV$-Gorenstein} objects, and associated categories
$\catKLtac{\capUV}$, $\catGL{\sfA}$, and $\stabcatGL{\sfA}$.

It is straightforward to verify that the functors described in the
next theorem are those induced by the embeddings of $\sfV$ and $\sfU$
into $\bcatD{b}{\sfV}$ and $\bcatD{b}{\sfU}$.

\begin{thm}
  \label{thm:BH-functor}
  Let $(\sfU,\sfV)$ be a complete cotorsion pair in $\sfA$.  There are
  fully faithful triangulated functors
  \begin{equation*}
    \dmapdef{\fuFPrj}{\stabcatGR{\sfA}}{\catDsgPrj{\sfV}} \qqand
    \dmapdef{\fuFInj}{\stabcatGL{\sfA}}{\catDsgInj{\sfU}} \:,
  \end{equation*}
  where $\fuFPrj$ sends a right $\sfU$-Gorenstein object to the hard
  truncation below at $0$ of a defining $\sfU$-totally acyclic
  complex, and $\fuFInj$ is defined similarly.
\end{thm}

\begin{prf*}
  We prove the assertion regarding $\fuFPrj$; the one for $\fuFInj$ is
  proved similarly. The functor
  $\catGR[\sfU]{\sfA} \to \catKRtac[\sfU]{\capUV}$ that maps a right
  $\sfU$-Gorenstein object to a defining right $\sfU$-totally acyclic
  complex induces a triangulated equivalence
  $\stabcatGR[\sfU]{\sfA} \to \catKRtac[\sfU]{\capUV}$; this was shown
  in \thmcite[3.8]{CET-20}. The arguments in the proof of
  \thmcite[3.1]{BJO-15} show that hard truncation below at $0$ yields
  a fully faithful triangulated functor
  $\catKRtac[\sfU]{\capUV} \to
  \bcatKeeac[\sfU]{+}{\capUV}/\bcatK{b}{\capUV}$.  Propositions
  \prpref[]{semi-exact} and \prpref[]{UV-exact} now yield the desired
  fully faithful triangulated functor $\fuFPrj$.
\end{prf*}

\thmref{BH-functor} facilitates the following definition; see also
\cite{BJO-15}.

\begin{dfn}
  \label{dfn:Gorenstein-defect}
  Let $(\sfU,\sfV)$ be a complete cotorsion pair in $\sfA$.  The
  Verdier quotient of the projective singularity category by the
  essential image of $\fuFPrj$,
  \begin{equation*}
    \catDdefPrj{\sfV} \colonequals \catDsgPrj{\sfV}/\im(\fuFPrj) \:,
  \end{equation*}
  is called the \emph{projective Gorenstein defect category of} $\sfV$
  or \emph{$\catPrj{\sfV}$-Gorenstein defect category}. Similarly, the
  \emph{injective Gorenstein defect category of} $\sfU$ or
  \emph{$\catInj{\sfU}$-Gorenstein defect category} is the Verdier
  quotient
  \begin{equation*}
    \catDdefInj{\sfU} \colonequals \catDsgInj{\sfU}/\im(\fuFInj) \:.
  \end{equation*}
\end{dfn}

It was shown in \lemcite[2.10]{CET-20} that the categories
$\catGR{\sfA}$ and $\catGL{\sfA}$ are closed under extensions. In the
context of a complete hereditary cotorsion pair, we now show that they
are closed under direct summands; our proof relies of work of Chen,
Liu, and Yang \cite{CLY-19,XYnWCh17}.  In the next two proofs we use
the notation $\mathcal{G}(\sfU)$ from \dfncite[3.1]{XYnWCh17} to
denote cycle subobjects of acyclic complexes of objects from $\sfU$
that are $\hom_\sfA(-,\capUV)$-exact.

\begin{lem}
  \label{lem:smdclosed}
  Let $(\sfU,\sfV)$ be a complete hereditary cotorsion pair in $\sfA$.
  The categories $\catGR{\sfA}$ and $\catGL{\sfA}$ are closed under
  direct summands.
\end{lem}

\begin{prf*}
  It follows from \dfncite[1.1]{CET-20} and \lemcite[3.2]{XYnWCh17}
  that the class $\catGR{\sfA}$ is the intersection of $\sfV$ and the
  class $\mathcal{G}(\sfU)$ defined in \dfncite[3.1]{XYnWCh17}. (We
  notice that \cite[3.1--3.3]{XYnWCh17} do not depend on the blanket
  assumption that the underlying abelian category is bicomplete.) It
  now follows from \prpcite[3.3]{CLY-19} that $\catGR{\sfA}$ is closed
  under direct summands. That $\catGL{\sfA}$ is closed under direct
  summands is proved similarly.
\end{prf*}

\begin{lem}
  \label{lem:sesclosed}
  Let $(\sfU,\sfV)$ be a complete hereditary cotorsion pair in $\sfA$.
  \begin{prt}
  \item Let $0\to V'\to V\to V''\to 0$ be an exact sequence in $\sfV$.
    \begin{rqm}
    \item If $V''$ is in $\catGR{\sfA}$, then $V'$ is in
      $\catGR{\sfA}$ if and only if $V$ is in $\catGR{\sfA}$.
    \item If $V'$ and $V$ are in $\catGR{\sfA}$, then
      $V'' \in \catGR{\sfA}$ if and only if
      $\Ext[\sfA]{1}{V''}{W} = 0$ holds for all $W\in{\sfU\cap\sfV}$.
    \end{rqm}
  \item Let $0\to U'\to U\to U''\to 0$ be an exact sequence in $\sfU$.
    \begin{rqm}
    \item If $U'$ is in $\catGL{\sfA}$, then $U''$ is in
      $\catGL{\sfA}$ if and only if $U$ is in $\catGL{\sfA}$.
    \item If $U$ and $U''$ are in $\catGL{\sfA}$, then $U'$ is in
      $\catGL{\sfA}$ if and only if $\Ext[\sfA]{1}{W}{U'} = 0$ holds
      for all $W\in{\sfU\cap\sfV}$.
    \end{rqm}
  \end{prt}
\end{lem}

\begin{prf*}
  As noticed in the proof of \lemref{smdclosed}, the class
  $\catGR{\sfA}$ is the intersection of $\sfV$ and the class
  $\mathcal{G}(\sfU)$ defined in \dfncite[3.1]{XYnWCh17}.  The
  assertions in part (a) now follow from \prpcite[3.3]{XYnWCh17}; the
  proof of part (b) is similar.
\end{prf*}

The next definition foreshadows notions of Gorenstein dimensions
relative to a cotorsion pair.

\begin{dfn}
  Let $(\sfU,\sfV)$ be a complete hereditary cotorsion pair in $\sfA$.
  \begin{prt}
  \item A complex $V$ in $\catDb{\sfV}$ is called
    \emph{$\catGR{\sfA}$-perfect} if $V$ is isomorphic in
    $\catDb{\sfV}$ to a complex $X$ with $X_n\in \catGR{\sfA}$ for all
    $n\in \ZZ$ and $X_n=0$ for $|n|\gg0$.
  \item A complex $U$ in $\catDb{\sfU}$ is called
    \emph{$\catGL{\sfA}$-perfect} if $U$ is isomorphic in
    $\catDb{\sfU}$ to a complex $Y$ with $Y_n\in \catGL{\sfA}$ for all
    $n\in \ZZ$ and $Y_n=0$ for $|n|\gg0$.
  \end{prt}
\end{dfn}

\begin{thm}
  \label{thm:LGor-dimension}
  Let $(\sfU,\sfV)$ be a complete hereditary cotorsion pair in $\sfA$.
  \begin{prt}
  \item For a complex $V \in \catDb{\sfV}$ the following conditions
    are equivalent:
    \begin{eqc}
    \item $V$ is $\catGR{\sfA}$-perfect.
    \item For every complex $X\in \bcatKeeac[\sfV]{+}{\capUV}$
      isomorphic to $V$ in $\catDb{\sfV}$ one has
      $\Cy[n]{X}\in \catGR{\sfA}$ for $n\gg0$.
    \item There exists a complex $X\in \bcatKeeac[\sfV]{+}{\capUV}$
      isomorphic to $V$ in $\catDb{\sfV}$ with
      $\Cy[n]{X}\in \catGR{\sfA}$ for $n\gg0$.
    \item For every complex $X\in \bcatKeeac[\sfV]{+}{\catGR{\sfA}}$
      isomorphic to $V$ in $\catDb{\sfV}$ one has
      $\Cy[n]{X}\in{\catGR{\sfA}}$ for $n\gg0$.
    \end{eqc}

  \item For a complex $U \in \catDb{\sfU}$ the following conditions
    are equivalent:
    \begin{eqc}
    \item $U$ is $\catGL{\sfA}$-perfect.
    \item For every complex $Y\in \bcatKeeac[\sfU]{$-$}{\capUV}$
      isomorphic to U in $\catDb{\sfU}$ one has
      $\Cy[n]{Y}\in \catGL{\sfA}$ for $n\ll0$.
    \item There exists a complex $Y\in \bcatKeeac[\sfU]{$-$}{\capUV}$
      isomorphic to $U$ in $\catDb{\sfU}$ with
      $\Cy[n]{Y}\in \catGL{\sfA}$ for $n\ll 0$.
    \item For every complex $Y\in \bcatKeeac[\sfU]{$-$}{\catGL{\sfA}}$
      isomorphic to $U$ in $\catDb{\sfU}$ one has
      $\Cy[n]{Y}\in{\catGL{\sfA}}$ for $n\ll 0$.
    \end{eqc}
  \end{prt}
\end{thm}

\begin{prf*}
  We prove (a); the proof of (b) is dual.

  \proofofimp{i}{ii} Let $X$ be a complex in
  $\bcatKeeac[\sfV]{+}{\capUV}$ with $V\cong X$ in $\catDb{\sfV}$. As
  $X$ is isomorphic in $\catK{\capUV}$ to a bounded below complex, and
  as $\capUV$ is closed under direct summands, the complex $X$ is
  isomorphic in $\catK{\capUV}$ to a soft truncation below, so we may
  assume that $X_n = 0$ holds for $n \ll 0$. By assumption, there is a
  complex $X'$ isomorphic to $V$ in $\catDb{\sfV}$ with
  $X'_n\in{\catGR{\sfA}}$ for all $n\in \ZZ$ and $X'_n=0$ for
  $|n|\gg0$. As $X\cong X'$ in $\catDb{\sfV}$, there exists a
  $\sfV$-quasi-isomorphism $\mapdef{\gra}{X}{X'}$; indeed, this
  follows from the dual statement to B\"uhler's
  \lemcite[3.3.3]{TBh11}, cf.~\prpref{UV-exact}.
  In particular, $\Cone{\gra}$ is a $\sfV$-acyclic complex. The
  category $\catGR{\sfA}$ is additive, so $\Cone{\gra}$ is a complex
  of right $\sfU$-Gorenstein objects and $\Cone{\gra}_n=0$ holds for
  $n\ll0$. It now follows from \lemref{sesclosed}(a,1) that
  $\Cy[n]{\Cone{\gra}}$ is right $\sfU$-Gorenstein for all $n\in
  \ZZ$. For $n\gg0$ one has $\Cy[n]{\Cone{\gra}} \cong \Cy[n]{X}$ and
  thus $\Cy[n]{X}$ belongs to $\catGR{\sfA}$.

  \proofofimp{ii}{iii} Trivial in view of \prpref{semi-exact}.

  \proofofimp{iii}{iv} Let $X$ be a complex in
  $\bcatKeeac[\sfV]{+}{\catGR{\sfA}}$ with $V\cong X$ in
  $\catDb{\sfV}$. By assumption, there exists a complex $X'$ in
  $\bcatKeeac[\sfV]{+}{\capUV}$ such that $X\cong X'$ in
  $\catDb{\sfV}$ and $\Cy[n]{X'}\in \catGR{\sfA}$ for $n\gg0$.  By
  \lemref{smdclosed} the category $\catGR{\sfA}$ is closed under
  direct summands, and so is $\capUV$. It follows that $X$ and $X'$
  are isomorphic to soft truncations below, so without loss of
  generality we assume that $X_n=0=X'_n$ holds for $n \ll 0$. As above
  there exists a $\sfV$-quasi-isomorphism $\mapdef{\gra}{X'}{X}$. Thus
  $\Cone{\gra}$ is a bounded below $\sfV$-acyclic complex of right
  $\sfU$-Gorenstein objects, so $\Cy[n]{\Cone{\gra}}$ is right
  $\sfU$-Gorenstein for every $n\in \ZZ$ by
  \lemref{sesclosed}(a,1). For $n \gg0$ the exact sequence
  $0\to X\to \Cone{\gra} \to \Shift{X'} \to 0$ now yields an exact
  sequence
  $0\to \Cy[n]{X} \to \Cy[n]{\Cone{\gra}}\to \Cy[n-1]{X'}\to 0$, and
  another application of \lemref{sesclosed}(a,1) yields
  $\Cy[n]{X} \in \catGR{\sfA}$.

  \proofofimp{iv}{i} In view of \prpref{semi-exact} there exists a
  complex $X$ in $\bcatKeeac[\sfV]{+}{\capUV}$ with $V\cong X$ in
  $\catDb{\sfV$}, and as above we can without loss of generality
  assume that $X_n=0$ holds for $n \ll 0$. By \eqclbl{iv} the cycle
  object $\Cy[n]{X}$ is right $\sfU$-Gorenstein for some $n \gg 0$, so
  the complex $0\to \Cy[n]{X} \to X_n\to \cdots$ is a bounded complex
  that is isomorphic to $V$ in $\catDb{\sfV}$ and whose objects belong
  to $\catGR{\sfA}$.
\end{prf*}

\begin{dfn}
  \label{dfn:RGLG}
  Let $(\sfU,\sfV)$ be a complete hereditary cotorsion pair in $\sfA$.
  We consider the following subcategories of $\catDb{\sfV}$ and
  $\catDb{\sfU}$:
  \begin{align*}
    \catfinR{\sfV} &\colonequals \{V \in \catDb{\sfV} \mid
                     \text{$V$ is $\catGR{\sfA}$-perfect}\}\quad\text{and}\\
    \catfinL{\sfU} &\colonequals\{U\in \catDb{\sfU} \mid
                     \text{$U$ is $\catGL{\sfA}$-perfect}\}\:.
  \end{align*}
\end{dfn}

\begin{prp}
  \label{prp:fGor_smallest_subcat}
  Let $(\sfU,\sfV)$ be a complete hereditary cotorsion pair in $\sfA$.
  \begin{prt}
  \item The category $\catfinR{\sfV}$ is thick in $\catDb{\sfV}$, and
    it is the smallest triangulated subcategory of $\catDb{\sfV}$ that
    contains the objects from $\catGR{\sfA}$.
  \item The category $\catfinL{\sfU}$ is thick in $\catDb{\sfU}$, and
    it is the smallest triangulated subcategory of $\catDb{\sfU}$ that
    contains the objects from $\catGL{\sfA}$.
  \end{prt}
\end{prp}

\begin{prf*}
  We prove part (a); the proof of (b) is similar.  Evidently
  $\catfinR{\sfV}$ is additive and closed under shifts and
  isomorphisms.

  Let $V'\to V \to V'' \to \Shift{V'}$ be a triangle in $\catDb{\sfV}$
  with $V',V \in \catfinR{\sfV}$. By \prpref{semi-exact}, there is a
  triangle $X' \xra{\gra} X \lra \Cone{\gra} \lra \Shift{X'}$ in
  $\bcatKeeac[\sfV]{+}{\capUV}$ and an isomorphism in $\catDb{\sfV}$
  of triangles:
  \begin{equation*}
    \xymatrix{V'\ar[r] \ar[d]^{\cong}
      & V\ar[r] \ar[d]^{\cong}
      & V''\ar[r] \ar[d]^{\cong}
      & \Shift{V'}\ar[d]^{\cong}\\
      X'\ar[r] & X\ar[r] & \Cone{\gra}\ar[r] & \Shift{X'} \:.}
  \end{equation*}
  Since the category $\catfinR{\sfV}$ is closed under isomorphisms and
  $V'$ and $V$ belong to $\catfinR{\sfV}$, so do $X$ and
  $X'$. Consequently, there is an exact sequence
  \begin{equation*}
    0 \lra X \lra \Cone{\gra} \lra \Shift{X'} \lra 0
  \end{equation*}
  of complexes in $\catC{\sfV}$.  The sequence
  $0\to \Cy[n]{X}\to \Cy[n]{\Cone{\gra}}\to \Cy[n]{\Shift{X'}}\to 0$
  is exact for all $n \gg 0$, and $\Cy[n]{X}$ and $\Cy[n]{\Shift{X'}}$
  belong to $\catGR{\sfA}$ by \thmref{LGor-dimension}. It now follows
  from \lemcite[2.10]{CET-20} that $\Cy[n]{\Cone{\gra}}$ belongs to
  $\catGR{\sfA}$ for $n\gg0$. Thus $\Cone{\gra}$ belongs to
  $\catfinR{\sfV}$, again by \thmref{LGor-dimension}. It follows that
  $\catfinR{\sfV}$ is a triangulated subcategory of $\catDb{\sfV}$.

  The triangulated subcategory $\catfinR{\sfV}$ holds the objects from
  $\catGR{\sfA}$. Let $\langle \catGR{\sfA} \rangle$ be the smallest
  triangulated subcategory of $\catDb{\sfV}$ that contains the objects
  from $\catGR{\sfA}$, one then has
  $\langle \catGR{\sfA} \rangle \subseteq \catfinR{\sfV}$. Let $V$ be
  an object in $\catfinR{\sfV}$. There exists a bounded complex $X$
  with $V \cong X$ in $\catDb{\sfV}$ such that each $X_n$ is in
  $\catGR{\sfA}$. Without loss of generality we may assume that one
  has $X_0\not=0$ and $X_n=0$ for $n<0$. Set
  $s=\sup\{n\mid X_n\not=0\}$; we proceed by induction on $s$. If
  $s=0$, clearly $X\in \langle\catGR{\sfA}\rangle$. For $s>0$, there
  is a triangle in $\catDb{\sfV}$
  \begin{equation*}
    \Tha{s-1}{X} \lra X \lra \Shift[s]{X_{s}} \lra \Shift{\Tha{s-1}{X}} \:.
  \end{equation*}
  By the induction hypothesis, both $\Tha{s-1}{X}$ and
  $\Shift[s]{X_s}$ belong to the triangulated category
  $\langle \catGR{\sfA}\rangle$, and hence so does $X$.

  It remains to show that $\catfinR{\sfV}$ is closed under direct
  summands. Let $V, V'\in\catDb{\sfV}$ and assume that $V\oplus V'$ is
  $\catGR{\sfA}$-perfect. There are $\sfV$-quasi-isomorphisms $X\to V$
  and $X'\to V'$ in $\catK{\sfV}$ with
  $X, X'\in\bcatKeeac[\sfV]{+}{\capUV}$, this follows from the dual to
  \lemcite[4.1]{BKl90}. It follows that $X\oplus X' \to V \oplus V'$
  is a $\sfV$-quasi-isomorphism in $\catK{\sfV}$. As $V\oplus V'$ is
  $\catGR{\sfA}$-perfect and
  $X\oplus X'\in\bcatKeeac[\sfV]{+}{\capUV}$, \thmref{LGor-dimension}
  implies that $\Cy[n]{X}\oplus \Cy[n]{X'}\cong \Cy[n]{X\oplus X'}$
  belongs to $\catGR{\sfA}$ for $n\gg0$. Thus $\Cy[n]{X}$ and
  $\Cy[n]{X'}$ belong to $\catGR{\sfA}$ for $n\gg0$ by
  \lemref{smdclosed}. Again from \thmref{LGor-dimension} it follows
  that $V$ and $V'$ belong to $\catfinR{\sfV}$.
\end{prf*}

Given \thmcite[2.11]{CET-20}, \thmref{LGor-dimension}, and
\prpref{fGor_smallest_subcat} one could obtain the first equivalence
in part (a) below from a result of Iyama and Yang
\corcite[2.2]{OImYDn20}, which elaborates on an example by Keller and
Vossieck \cite{BKlDVs87}.

\begin{thm}
  \label{thm:Gorenstein-defect}
  Let $(\sfU,\sfV)$ be a complete hereditary cotorsion pair in $\sfA$.
  \begin{prt}
  \item There are triangulated equivalences
    \begin{equation*}
      \stabcatGR{\sfA} \dqis \frac{\catfinR{\sfV}}{\catKb{\capUV}} \qqand
      \catDdefPrj{\sfV} \dqis \frac{\catDb{\sfV}}{\catfinR{\sfV}} \:.
    \end{equation*}
  \item There are triangulated equivalences
    \begin{equation*}
      \stabcatGL{\sfA} \dqis \frac{\catfinL{\sfU}}{\catKb{\capUV}} \qqand
      \catDdefInj{\sfU} \dqis \frac{\catDb{\sfU}}{\catfinL{\sfU}} \:.
    \end{equation*}
  \end{prt}
\end{thm}

\begin{prf*}
  We prove part (a); the equivalences in part (b) have similar
  proofs. By \prpref{fGor_smallest_subcat} it suffices to show that
  the quotient $\catfinR{\sfV}/\catKb{\capUV}$ is the essential image
  of the functor
  $\fuFPrj:\stabcatGR{\sfA}\to \catDb{\sfV}/\catKb{\capUV}$. Let
  $V\in{\catfinR{\sfV}}$, by \prpref{semi-exact} there is an
  isomorphism $V\cong X$ in $\catDb{\sfV}$ for some
  $X\in \bcatKeeac[\sfV]{+}{\capUV}$. By \thmref{LGor-dimension} there
  exists an integer $n$ such that $\Cy[n]{X}\in \catGR{\sfA}$ and
  $\Tsb{n}{X}$ is $\sfV$-acyclic; indeed, the cycle objects are right
  $\sfU$-Gorenstein by \lemref{sesclosed}(a,1).  Let $\widetilde{T}$
  be a right $\sfU$-totally acyclic complex with
  $\Cy[0]{\widetilde{T}}=\Cy[n]{X}$. Splicing together
  $\Shift[n]{\Tha{0}{\widetilde{T}}}$ and $\Thb{n+1}{X}$ we obtain a
  right $\sfU$-totally acyclic complex $T$ with
  $\Thb{n+1}{T} = \Thb{n+1}{X}$. By definition one has
  $\fuFPrj(\Cy[0]{T})=\Thb{0}{T}$, and there are isomorphisms
  $\Thb{0}{T}\cong \Thb{n+1}{T}\cong X$ in $\catDsgPrj{\sfV}$.

  Since the projective objects in $\sfV$ are those in $\capUV$, see
  \prpref{UV-exact}, the equivalence
  $\catDdefPrj{\sfV}\simeq\catDb{\sfV}/\catfinR{\sfV}$ follows from
  \dfnref{Dsg}(a) and \dfnref{Gorenstein-defect}.
\end{prf*}

\begin{cor}
  \label{cor:A.6}
  Let $(\sfU,\sfV)$ be a complete hereditary cotorsion pair in $\sfA$.
  \begin{prt}
  \item The following conditions are equivalent:
    \begin{eqc}
    \item $\mapdef{\fuFPrj}{\stabcatGR{\sfA}}{\catDsgPrj{\sfV}}$ is a
      triangulated equivalence.
    \item One has $\catDdefPrj{\sfV} =0$.
    \item Every object in $\catDb{\sfV}$ is $\catGR{\sfA}$-perfect.
    \item Every object in $\sfV$ is $\catGR{\sfA}$-perfect.
    \end{eqc}
  \item The following conditions are equivalent:
    \begin{eqc}
    \item $\mapdef{\fuFInj}{\stabcatGL{\sfA}}{\catDsgInj{\sfU}}$ is a
      triangulated equivalence.
    \item One has $\catDdefInj{\sfU} =0$.
    \item Every object in $\catDb{\sfU}$ is $\catGL{\sfA}$-perfect.
    \item Every object in $\sfU$ is $\catGL{\sfA}$-perfect.
    \end{eqc}
    \noindent
  \end{prt}
\end{cor}

\begin{prf*}
  We prove part (a); the proof of part (b) is similar.  The
  equivalence of the first three conditions follows immediately from
  \dfnref{Gorenstein-defect} and \thmref{Gorenstein-defect}. The
  equivalence of \eqclbl{iii} and \eqclbl{iv} follows from
  \prpref{fGor_smallest_subcat}, as the smallest triangulated
  subcategory of $\catDb{\sfV}$ that contains $\sfV$ is
  $\catDb{\sfV}$.
\end{prf*}

\begin{rmk}
  \label{rmk:DSK}
  Let $\sfUV$ be a cotorsion pair in $\sfA$. When one combines the
  previous corollary with the equivalence from \thmcite[3.8]{CET-20},
  it transpires that if every object in $\sfV$ is
  $\catGR{\sfA}$-perfect, then there are triangulated equivalences
  \begin{equation*}
    \catDsgPrj{\sfV} \dqis \stabcatGR{\sfA} \dqis \catKRtac{\capUV} \:.
  \end{equation*}
  Similarly, if every object in $\sfU$ is $\catGL{\sfA}$-perfect, then
  there are triangulated equivalences
  \begin{equation*}
    \catDsgInj{\sfU} \dqis \stabcatGL{\sfA} \dqis \catKLtac{\capUV} \:.
  \end{equation*}
\end{rmk}

In \cite{CET-20} right and left totally acyclic complexes and right
and left Gorenstein objects are defined for any subcategory of $\sfA$.
Given a cotorsion pair $\sfUV$ in $\sfA$, the subcategory $\capUV$ is
self-orthogonal, and for such a category the right/left distinctions
disappear, see \dfncite[1.6 and 2.1]{CET-20}, and one speaks of
$\capUV$-totally acyclic complexes and $\capUV$-Gorenstein objects.
The homotopy category of $\capUV$-totally acyclic complexes is denoted
$\catKtac{\capUV}$ and the category of $\capUV$-Gorenstein objects is
denoted $\catG[\capUV]{\sfA}$.

\begin{prp}
  \label{prp:intersect}
  Let $(\sfU,\sfV)$ be a complete hereditary cotorsion pair in $\sfA$.
  \begin{prt}
  \item The following conditions are equivalent:
    \begin{eqc}
    \item $\catG[\sfU\cap\sfV]{\sfA}\subseteq \sfV$.
    \item $\catG[\sfU\cap\sfV]{\sfA}=\catGR{\sfA}$.
    \end{eqc}
    Moreover, when these conditions hold one has
    $\catGL{\sfA}=\capUV$.
  \item The following conditions are equivalent:
    \begin{eqc}
    \item $\catG[\sfU\cap\sfV]{\sfA}\subseteq \sfU$.
    \item $\catG[\sfU\cap\sfV]{\sfA}=\catGL{\sfA}$.
    \end{eqc}
    Moreover, when these conditions hold one has
    $\catGR{\sfA}=\capUV$.
  \end{prt}
\end{prp}

\begin{prf*}
  We prove (a); the proof of (b) is similar.

  \proofofimp{ii}{i} The containment $\catGR{\sfA}\subseteq \sfV$
  holds by definition.

  \proofofimp{i}{ii} The containment
  $\catG[\sfU\cap\sfV]{\sfA}\supseteq \catGR{\sfA}$ holds by
  \remcite[1.8]{CET-20}. If one has
  $\catG[\sfU\cap\sfV]{\sfA}\subseteq \sfV$, then every
  $\capUV$-totally acyclic complex is right $\sfU$-totally acyclic,
  and hence \eqclbl{ii} follows.

  To see that the moreover statement holds under these conditions,
  recall from \dfncite[1.1 and Exa.~1.2]{CET-20} that there are
  containments $\capUV \subseteq \catGL{\sfA}\subseteq \sfU$. Further,
  $\catGL{\sfA}\subseteq \catG[\sfU\cap\sfV]{\sfA}$ holds by
  \remcite[1.8]{CET-20}, so in view of $(i)$ one now has
  $\catGL{\sfA} \subseteq \sfV$ and hence $\catGL{\sfA} = \capUV$.
\end{prf*}

\begin{rmk}
  \label{rmk:intersect}
  Let $\sfUV$ be a cotorsion pair in $\sfA$.  For the self-orthogonal
  class $\capUV$, the equivalences in \rmkref{DSK} simplify as
  follows:

  If every object in $\sfV$ is $\catGR{\sfA}$-perfect and
  $\catG[\sfU\cap\sfV]{\sfA}\subseteq \sfV$, then there are
  triangulated equivalences
  \begin{equation*}
    \catDsgPrj{\sfV} \dqis \stabcatG[\capUV]{\sfA} \dqis \catKtac{\capUV} \:.
  \end{equation*}

  Similarly, if every object in $\sfU$ is $\catGL{\sfA}$-perfect and
  $\catG[\sfU\cap\sfV]{\sfA}\subseteq \sfU$, then there are
  triangulated equivalences
  \begin{equation*}
    \catDsgInj{\sfU} \dqis \stabcatG[\capUV]{\sfA} \dqis \catKtac{\capUV} \:.
  \end{equation*}
\end{rmk}

\begin{exa}
  Let $A$ be a ring and consider the cotorsion pair
  $(\catPrj{A},\catMod{A})$. A $\catPrj{A}$-Gorenstein module is a
  Gorenstein projective module in the classic sense, see
  \exacite[2.5]{CET-20}. If every $A$-module has finite Gorenstein
  projective dimension, then the categories in the first display in
  \rmkref{intersect} are ``big'' versions of the equivalent categories
  $(\diamond)$ from the introduction.
\end{exa}

\section{Gorenstein dimensions for everyone}
\label{sec:gdim}

\noindent
Assume that $\sfA$ is Grothendieck and $(\sfU, \sfV)$ is a complete
hereditary cotorsion pair in $\sfA$.  Under this extra, but not too
restrictive, assumption on $\sfA$ one can extend the notions of
$\catGR{\sfA}$- and $\catGL{\sfA}$-perfection to the derived category
of $\sfA$: To every complex $M$ of objects from $\sfA$ we assign two
numbers, the $\catGR{\sfA}$-projective dimension and the
$\catGL{\sfA}$-injective dimension, and for objects in $\catDb{\sfV}$
and $\catDb{\sfU}$ these dimensions are finite if and only if the
objects are $\catGR{\sfA}$- and $\catGL{\sfA}$-perfect, respectively.

The category $\catC{\sfA}$ of chain complexes is also Grothendieck,
and it follows from work of Gillespie \prpcite[3.6]{JGl04} and \JSt\
\prpcite[7.14]{JSt13b} that $(\sfU,\sfV)$ induces two complete
hereditary cotorsion pairs in $\catC{\sfA}$:
\begin{equation}
  \label{eq:pairs}
  (\textsf{semi-}\sfU,\sfV\textsf{-ac}) \qqand
  (\sfU\textsf{-ac},\textsf{semi-}\sfV) \:.
\end{equation}
The complexes in \textsf{$\sfV$-ac} are acyclic complexes of objects
in $\sfV$ with cycle objects in $\sfV$, i.e.\ the $\sfV$-acyclic
complexes in $\catK{\sfV}$. The semi-$\sfU$ complexes, i.e.\ the
objects in \textsf{semi-$\sfU$}, are complexes $U$ of objects in
$\sfU$ with the property that $\Hom[\sfA]{U}{V}$ is acyclic for every
$\sfV$-acyclic complex $V$. The classes \textsf{$\sfU$-ac} and
\textsf{semi-$\sfV$} are defined similarly.  We call a complex in
$\textsf{semi-}\sfU \cap \textsf{semi-}\sfV$ a \emph{\semiUV\
  complex.}

Completeness of the cotorsion pairs in \eqref{pairs} yields:

\begin{fct}
  \label{fct:complete}
  For every $\sfA$-complex $M$ there are exact sequences of
  $\sfA$-complexes,
  \begin{equation*}
    0 \lra V' \lra U \xra{\pi} M \lra 0 \qqand
    0 \lra M \xra{\iota} V \lra U' \lra 0 \:,
  \end{equation*}
  where $U$ is $\textsf{semi-}\sfU$, $V$ is $\textsf{semi-}\sfV$, $V'$
  is $\sfV$-acyclic, $U'$ is $\sfU$-acyclic, and the homomorphisms
  $\pi$ and $\iota$ are quasi-isomorphisms. See \prpcite[3.6]{JGl04}
  and \prpcite[7.14]{JSt13b}.
\end{fct}

\begin{exa}
  \label{exa:bd}
  A complex $U$ of objects in $\sfU$ with $\Co[n]{U}\in \sfU$ for
  $n\ll0$ is a semi-$\sfU$ complex and a complex $V$ of objects in
  $\sfV$ with $\Cy[n]{V}\in \sfV$ for $n\gg0$ is a semi-$\sfV$
  complex. This follows from \cite[Props. A.1 and
  A.3]{LWCPTh19}\footnote{The proofs of \cite[Props. A.1 and
    A.3]{LWCPTh19} apply to complexes of objects in any abelian
    category.}.

  In particular, a bounded below complex of objects in $\sfU$ is a
  semi-$\sfU$ complex and a bounded above complex of objects in $\sfV$
  is a semi-$\sfV$ complex; cf. \lemcite[3.4]{JGl04}.
\end{exa}

\begin{dfn}
  \label{dfn:dgU-V-replacement}
  Let $M$ be an $\sfA$-complex. A \emph{semi-$\sfU$-$\sfV$ replacement
    of $M$} is a semi-$\sfU$-$\sfV$ complex that is isomorphic to $M$
  in the derived category $\catD{\sfA}$.
\end{dfn}

Some technical results about the cotorsion pairs \eqref{pairs} and
\semiUV\ complexes have been relegated to an appendix. The first step
towards defining the $\catGR{\sfA}$-projective and
$\catGL{\sfA}$-injective dimensions is to notice that every
$\sfA$-complex has a \semiUV\ replacement.

\begin{rmk}
  \label{rmk:replace}
  Per \lemref{approx}(a) every semi-$\sfV$ complex $V$ has a \semiUV\
  replacement, as there is even a surjective quasi-isomorphism
  $W \qra V$ with $W$ \semiUV. Similarly, by \lemref{approx}(b), every
  semi-$\sfU$ complex $U$ has a \semiUV\ replacement as there is even
  an injective quasi-isomorphism $U \qra W'$ with $W'$ \semiUV. These
  combine to show that every $\sfA$-complex $M$ has a \semiUV\
  replacement:

  Let $M$ be an $\sfA$-complex. By \fctref{complete} there is a
  semi-$\sfV$ complex $V$ and a quasi-isomorphism $M \qra V$, and
  \lemref{approx}(a) yields a \semiUV\ complex $W$ and a
  quasi-isomorphism $W \qra V$. Thus $W$ is a \semiUV\ replacement of
  $M$.

  One could also start with a quasi-isomorphism $U \qra M$ as in
  \lemref{approx}(a) and get a quasi-isomorphism $U \qra W'$ from
  \lemref{approx}(b).
\end{rmk}

The classic notion of Gorenstein projective dimension of modules and
complexes over a ring $A$ is a special case of the
$\catGR{\sfA}$-projective dimension defined below: Apply the
definition with $\sfA = \catMod{A}$ and
$(\sfU,\sfV) = (\catPrj{A},\catMod{A})$; see also
\exacite[1.7]{CET-20} and \rmkcite[4.8]{CELTWY-21}.

\begin{dfn}
  \label{dfn:Gfcd}
  Let $\sfUV$ be a complete hereditary cotorsion pair in $\sfA$ and
  $M$ an $\sfA$-complex.
  \begin{prt}
  \item The \emph{$\catGR{\sfA}$-projective dimension of $M$} is given
    by
    \begin{align*}
      \catGR{\sfA}&\operatorname{-pd} M \deq \\
                  & \inf%
                    \left\{n\in\ZZ
                    \:\left|\:
                    \begin{gathered}
                      \text{There is a semi-$\sfU$-$\sfV$ replacement
                        $W$ of $M$
                        with}\\[-3pt]
                      \text{$\H[i]{W}=0$ for all $i > n$ and
                        $\textrm{C}_{n}(W)$ in $\catGR{\sfA}$}
                    \end{gathered}
      \right.
      \right\}\hspace{-2pc}
    \end{align*}
    with $\inf\varnothing = \infty$ by convention.

  \item The \emph{$\catGL{\sfA}$-injective dimension of $M$} is given
    by
    \begin{align*}
      \catGL{\sfA}&\operatorname{-id} M \deq \\
                  & \inf%
                    \left\{n\in\ZZ
                    \:\left|\:
                    \begin{gathered}
                      \text{There is a semi-$\sfU$-$\sfV$ replacement $W$ of $M$ with}\\[-3pt]
                      \text{$\H[i]{W}=0$ for all $i < -n$ and
                        $\Cy[-n]{W}$ in $\catGL{\sfA}$}
                    \end{gathered}
      \right.
      \right\}\hspace{-2pc}
    \end{align*}
    with $\inf\varnothing = \infty$ by convention.
  \end{prt}
\end{dfn}

Note that
$\catGR{\sfA}\operatorname{-pd}M = -\infty =
\catGL{\sfA}\operatorname{-id}M$ holds if $M$ is acyclic.  We continue
with a caveat that compares to \rmkcite[5.9]{CELTWY-21}: An object of
$\catGR{\sfA}$-projective dimension $0$ need not belong to
$\catGR{\sfA}$; similarly for $\catGL{\sfA}$.

\begin{rmk}
  For every non-zero object $U \in \sfU$ one has
  $\catGR{\sfA}\operatorname{-pd}U = 0$. Indeed, the inequality
  ``$\ge$'' holds as $\H[0]{U} \ne 0$ and the opposite inequality
  holds as one can construct a \semiUV\ replacement of $U$
  concentrated in non-positive degrees: Completeness of $(\sfU,\sfV$)
  yields a coresolution of $U$ by objects in $\capUV$ that is
  concentrated in non-positive degrees and whose cycle objects are in
  $\sfU$; it is a semi-$\sfU$-$\sfV$ complex by \exaref{bd}.

  Similarly, $\catGL{\sfA}\operatorname{-id}V = 0$ holds for every
  non-zero object $V \in \sfV$.
\end{rmk}

\begin{lem}
  \label{lem:any-Co}
  Let $\sfUV$ be a complete hereditary cotorsion pair in $\sfA$ and
  $M$ an $\sfA$-complex with a \semiUV\ replacement $W$.
  \begin{prt}
  \item For every integer $n$ with
    $\catGR{\sfA}\operatorname{-pd}M \leq n$ one has
    $\Co[n]{W} \in \catGR{\sfA}$.
  \item For every integer $n$ with
    $\catGL{\sfA}\operatorname{-id}M \leq n$ one has
    $\Cy[-n]{W} \in \catGL{\sfA}$.
  \end{prt}
\end{lem}

\begin{prf*}
  We prove part (a); the proof of part (b) is similar. One can assume
  that $\catGR{\sfA}\operatorname{-pd}M = g$ holds for some integer
  $g$. By assumption there exists a \semiUV\ replacement $W'$ of $M$
  with $\Co[g]{W'}$ in $\catGR{\sfA}$ and $\H[n]{W'}=0$ for $n>g$. By
  induction it follows from \lemref{sesclosed}(a) that $\Co[n]{W'}$
  belongs to $\catGR{\sfA}$ for every $n \ge g$. Let $W$ be any
  \semiUV\ replacement of $M$. It follows from \prpref{schanuel} and
  \lemref{smdclosed} that $\Co[n]{W}$ belongs to $\catGR{\sfA}$ for
  every $n \ge g$.
\end{prf*}

\begin{thm}
  \label{thm:perfect}
  Let $\sfUV$ be a complete hereditary cotorsion pair in $\sfA$.
  \begin{prt}
  \item An object $V \in \catDb{\sfV}$ is $\catGR{\sfA}$-perfect if
    and only if it has finite $\catGR{\sfA}$-projective dimension.
  \item An object $U \in \catDb{\sfU}$ is $\catGL{\sfA}$-perfect if
    and only if it has finite $\catGL{\sfA}$-injective dimension.
  \end{prt}
\end{thm}

\begin{prf*}
  We prove part (a); the proof of part (b) is similar. Without loss of
  generality, let $V$ be a bounded complex of objects in $\sfV$. There
  exists by \thmref{bdres}(a) a bounded below \semiUV\ replacement $W$
  of $V$.  The assertion now follows from \thmref{LGor-dimension} and
  \lemref{any-Co} as one has $\Cy[n]{W} \is \Co[n+1]{W}$ for
  $n > \sup\{i\in \ZZ \mid \H[i]{V}\not=0\}$.
\end{prf*}

\begin{rmk}
  \label{rmk:vanishing}
  Let $\sfUV$ be a complete hereditary cotorsion pair in $\sfA$.
  Standard arguments, see \thmcite[4.5]{CELTWY-21}, based on
  \lemref{sesclosed}(a,2) and \lemref[]{sesclosed}(b,2) show that the
  $\catGR{\sfA}$-projective and $\catGL{\sfA}$-injective dimensions of
  an $\sfA$-complex $M$ when finite can be detected in terms of
  vanishing of cohomology:
  \begin{align*}
    \catGR{\sfA}\operatorname{-pd} M 
    &= \sup\{ m\in \ZZ \mid \Ext[\sfA]{m}{M}{W} \ne 0
      \text{ for some } W \in \capUV \} \quad\text{and}\\
    \catGL{\sfA}\operatorname{-id} M 
    &= \sup\{ m\in \ZZ \mid \Ext[\sfA]{m}{W}{M} \ne 0
      \text{ for some } W \in \capUV \} \:. \\
  \end{align*}
\end{rmk}

\section{Gorenstein rings and schemes}
\label{sec:appl1}

\noindent
Our main interest here is schemes, but we warm up with rings.

Let $A$ be a ring. The cotorsion pair $(\catFlat{A},\catCot{A})$ in
the Grothendieck category $\catMod{A}$ is complete and hereditary and
hence gives rise to the $\catGR[\sfFlat]{A}$-projective dimension. By
\prpcite[4.2]{CET-20} this cotorsion pair satisfies the conditions in
\prpref{intersect}(a), and the $\catGR[\sfFlat]{A}$-projective
dimension was studied in \cite{CELTWY-21} under the name ``Gorenstein
flat-cotorsion dimension;'' see also \rmkcite[4.8]{CELTWY-21}.  In
this section we keep with that terminology as it more descriptive, and
we write $\catFlatCot{A}$ for the intersection
$\catFlat{A} \cap \catCot{A}$.

If every $A$-module has finite Gorenstein flat-cotorsion dimension,
then, as noted in \rmkref{intersect}, there are triangulated
equivalences
\begin{equation*}
  \catDsgPrj{\catCot{A}} \dqis \stabcatG[\sfFlatCot]{A}
  \dqis \catKtac{\catFlatCot{A}} \:.
\end{equation*}
We proceed to identify a class of rings over which these equivalences
are realized.

Recall that a commutative noetherian ring is called
\emph{Iwanaga--Gorenstein} if it has finite self-injective
dimension. The next theorem collects characterizations of such rings
in terms of (Gorenstein) flat--cotorsion theory that closely mirror
classic characterizations in terms of (Gorenstein) projectivity. In
particular, $(ii)$ is the analogue of Buchweitz's equivalence between
the stable category of finitely generated Gorenstein projective
modules and the singularity category \thmcite[4.4.1]{ROB86}. Condition
$(iii)$ is analogous to \thmcite[3.6]{BJO-15} and compares to
\rmkcite[5.6]{SInHKr06}. Conditions $(iv)$ and $(v)$ compare to
Auslander and Bridger's \thmcite[4.20]{MAsMBr69}.

\begin{thm}
  \label{thm:Iwanaga_Gor_ring}
  Let $A$ be commutative noetherian; the next conditions are
  equivalent.
  \begin{eqc}
  \item $A$ is Iwanaga--Gorenstein.
  \item The functor
    $\fuFPrj\colon \stabcatG[\sfFlatCot]{A} \lra
    \catDsgPrj{\catCot{A}}$ yields a triangulated equivalence.
  \item One has $\catDdefPrj{\catCot{A}} \deq 0$.
  \item Every cotorsion $A$-module has finite Gorenstein
    flat-cotorsion dimension.
  \item Every $A$-module has finite Gorenstein flat-cotorsion
    dimension.
  \item Every $A$-complex with bounded above homology has finite
    Gorenstein flat-cotorsion dimension.
  \end{eqc}
\end{thm}

\begin{prf*}
  Conditions $(ii)$, $(iii)$, and $(iv)$ are equivalent by
  \corref{A.6} and \thmref{perfect}. Conditions $(i)$, $(v)$, and
  $(vi)$ are equivalent by \corcite[5.10]{CELTWY-21}. Evidently $(v)$
  implies $(iv)$. For the converse let $M$ be an $A$-module and
  consider an exact sequence $0 \to M \to C \to F \to 0$ where $C$ is
  cotorsion and $F$ is flat. By \corcite[5.8]{CELTWY-21} an $A$-module
  has finite Gorenstein flat-cotorsion dimension if and only if it has
  finite Gorenstein flat dimension. It now follows from Holm
  \thmcite[3.15]{HHl04a} that also $M$ has finite Gorenstein
  flat-cotorsion dimension.
\end{prf*}

Part of the allure of the category $\catDsgPrj{\catCot{A}}$ is that it
does not rely on projective modules and thus persists in the
non-affine setting.  Let $X$ denote a scheme with structure sheaf
$\calO_X$. By a sheaf on $X$ we mean a quasi-coherent sheaf, and
$\Qcoh{X}$ denotes the category of such sheaves. We say that $X$ is
\emph{semi-separated} if it has an open affine covering with the
property that the intersection of any two open affine sets is affine.

Let $X$ be semi-separated noetherian. The flat sheaves and cotorsion
sheaves on $X$ form a complete hereditary cotorsion pair
$(\catFlat{X},\catCot{X})$ in the Grothendieck category $\Qcoh{X}$;
see \rmkcite[2.4]{CET}. It follows from \thmcite[3.3]{CET} that this
pair satisfies the equivalent conditions in \prpref{intersect}(a), so
as above we refer to the $\catGR[\sfFlat]{X}$-projective dimension as
the \emph{Gorenstein flat-cotorsion dimension.}

If every cotorsion sheaf on $X$ has finite Gorenstein flat-cotorsion
dimension then, as noted in \rmkref{intersect}, there are triangulated
equivalences:
\begin{equation*}
  \catDsgPrj{\catCot{X}}
  \dqis \stabcatG[\sfFlatCot]{X} \dqis  \catKtac{\catFlatCot{X}} \:.
\end{equation*}
We proceed to identify a class of schemes over which these
equivalences are realized. Recall that a semi-separated noetherian
scheme $X$ is \emph{Gorenstein} provided that ${\mathcal O}_{X,x}$ is
a Gorenstein ring for every $x\in X$. If, in addition, the scheme has
finite Krull dimension, then this is equivalent to saying that
${\mathcal O}_{X}(U)$ is a Gorenstein ring for each open affine set
$U$ in some, equivalently every, open affine covering of $X$.

\begin{thm}
  \label{thm:mainX}
  Let $X$ be a semi-separated noetherian scheme of finite Krull
  dimension. The following conditions are equivalent.
  \begin{eqc}
  \item $X$ is Gorenstein.
  \item The functor
    \begin{equation*}
      \fuFPrj\colon \stabcatG[\sfFlatCot]{X} \lra \catDsgPrj{\catCot{X}}
    \end{equation*}
    yields a triangulated equivalence.
  \item One has $\catDdefPrj{\catCot{X}}=0$.
  \item Every cotorsion sheaf on $X$ has finite Gorenstein
    flat-cotorsion dimension.
  \item Every sheaf on $X$ has finite Gorenstein flat-cotorsion
    dimension.
  \item Every complex of sheaves on $X$ with bounded above homology
    has finite Gorenstein flat-cotorsion dimension.
  \end{eqc}
\end{thm}

\begin{prf*}
  Conditions $(ii)$, $(iii)$, and $(iv)$ are equivalent by
  \corref{A.6} and \thmref{perfect}. Evidently, $(vi)$ implies $(v)$
  implies $(iv)$, so it suffices to show that $(i)$ implies $(vi)$ and
  $(iv)$ implies $(i)$.  Denote by $d$ the Krull dimension of $X$.

  \proofofimp{i}{vi} Let $\M$ be a complex of sheaves on $X$ with
  bounded above homology and assume without loss of generality that
  $\H[n]{\M} = 0$ holds for $n > 0$. Let $\F$ be a semi-flat-cotorsion
  replacement of $\M$.  Notice that the truncated complex
  $\Thb{0}{\F}$ is a semi-flat-cotorsion replacement of the cokernel
  $\C = \Co[0]{\F}$; set $\G = \Cy[d-1]{\F}$ and consider the exact
  sequence
  \begin{equation*}
    0 \lra \G \lra \F_{d-1} \lra \cdots \lra \F_0 \lra \C \lra 0 \:.
  \end{equation*}
  The truncated complex $\Thb{d}{\F}$ yields a left resolution of $\G$
  by cotorsion sheaves; splicing this with a coresolution by injective
  sheaves, one obtains $\G$ as a cycle in an acyclic complex of
  cotorsion sheaves, whence $\G$ is cotorsion by
  \thmcite[3.3]{CET}. Per \thmcite[4.3]{CET} it now suffices to show
  that $\G$ is Gorenstein flat in the sense of \dfncite[1.2]{CET}. Let
  $\mathcal{U}$ be an open affine covering of $X$. For every open set
  $U\in \mathcal{U}$ there is an exact sequence of
  ${\mathcal O}_X(U)$-modules,
  \begin{equation*}
    0 \lra \G(U) \lra \F_{d-1}(U) \lra \cdots \lra \F_0(U) \lra
    \C(U) \lra 0 \:,
  \end{equation*}
  with $\F_n(U)$ a flat ${\mathcal O}_X(U)$-module for
  $0\le n \le d-1$. As ${\mathcal O}_X(U)$ is Gorenstein of Krull
  dimension at most $d$, the module $\G(U)$ is Gorenstein flat; see
  e.g.\ \thmcite[12.3.1]{rha}. From \thmcite[1.6]{CET} it follows that
  the sheaf $\G$ is Gorenstein flat.

  \proofofimp{iv}{i} Let $\calU$ be an open affine covering of $X$ and
  fix an open set $U\in\calU$. The goal is to prove that the ring
  $\calO_X(U)$ is Gorenstein. As it has finite Krull dimension, it
  suffices to show that every $\calO_X(U)$-module has finite
  Gorenstein flat dimension; see e.g.\ \thmcite[12.3.1]{rha}. Let $M$
  be an $\calO_X(U)$-module; denote by $\widetilde{M}$ the
  corresponding sheaf on $U$ and by $(i_U)_*(\widetilde{M})$ the sheaf
  on $X$ induced by the embedding $i_U\colon U \to X$. As
  $(\catFlat{X},\catCot{X})$ is a complete hereditary cotorsion pair
  in $\Qcoh{X}$, see \rmkcite[2.4]{CET}, there is a flat resolution
  $\F$ of $(i_U)_*(\widetilde{M})$ over $X$ constructed by taking
  special flat precovers. In particular, for $n\ge 1$ the sheaf
  $\Co[n]{\F}$ is cotorsion and $\F_n$ is flat-cotorsion.  By
  assumption $\Co[1]{\F}$ has finite Gorenstein flat-cotorsion
  dimension, which means that $\Co[n]{\F}$ is Gorenstein
  flat-cotorsion for some $n\ge 1$. In particular, $\Co[n]{\F}$ is by
  \thmcite[4.3]{CET} Gorenstein flat. Thus one gets an exact sequence
  of $\calO_X(U)$-modules,
  \begin{equation*}
    0 \lra \Co[n]{\F}(U) \lra \F_{n-1}(U) \lra \cdots \lra \F_0(U)
    \lra M \lra 0 \:,
  \end{equation*}
  which shows that $M$ has finite Gorenstein flat dimension,
  cf.~\thmcite[1.6]{CET}.
\end{prf*}

\begin{rmk}
  \label{rmk:gorX}
  Let $X$ be a Gorenstein scheme of finite Krull dimension. It follows
  from \corcite[4.6]{CET}, \thmcite[4.27]{DMfSSl11}, and Murfet's
  \thmcite[7.9]{DMf-phd} that the stable category
  $\stabcatG[\sfFlatCot]{X}$ is compactly generated. Further, the
  opposite category of Orlov's singularity category
  \dfncite[1.8]{DOO04} is equivalent, up to direct factors, to the
  subcategory $\stabcatG[\sfFlatCot]{X}^{\mathrm{c}}$ of compact
  objects. Thus, it follows from \thmref{mainX} that the singularity
  category $\catDsgPrj{\catCot{X}}$ is compactly generated and that
  the opposite of Orlov's singularity category, up to direct factors,
  is equivalent to $\catDsgPrj{\catCot{X}}^{\mathrm{c}}$.
\end{rmk}

Finally we justify that $\catDsgPrj{\catCot{X}}$ is a non-affine
avatar of $\catDsgPrj{A}$:

\begin{thm}
  Let $A$ be a commutative Gorenstein ring of finite Krull dimension
  and set $X=\text{Spec}(A)$. There are equivalences of triangulated
  categories
  \begin{gather*}
    \catDsgPrj{A} \dqis \stabcatG[\sfPrj]{A} \dqis
    \catKtac{\catPrj{A}}
    \\[-.7pc]
    \begin{rotate}{-90}$\qis$
    \end{rotate}\\[.15pc]
    \catDsgPrj{\catCot{X}} \dqis \stabcatG[\sfFlatCot]{X} \dqis
    \catKtac{\catFlatCot{X}} \:.
  \end{gather*}
\end{thm}

\begin{prf*}
  A commutative noetherian ring of finite Krull dimension is
  Gorenstein if and only if it is Iwanaga--Gorenstein.  The horizontal
  equivalences thus come from the two theorems above combined with
  \rmkref{intersect}. The vertical equivalence follows from
  \corcite[5.9]{CET-20}.
\end{prf*}

\section{Finitistic dimension}
\label{sec:appl2}

\noindent
As discussed in \secref{exact}, vanishing of singularity categories
can capture finiteness of global dimensions. In this section we show
how vanishing of a defect category captures finiteness of a finitistic
dimension.

In this section $A$ is a left noetherian ring.  The \emph{little
  finitistic dimension} of~$A$, $\findim{A}$, is the supremum of
projective dimensions of finitely generated $A$-modules of finite
projective dimension. It is conjectured that $\findim{A}<\infty$ holds
for Artin algebras $A$; see Bass \cite{HBs60} and
\cite[Conjectures]{rta}. Originally raised as a question by Rosenberg
and Zelinsky, this is now the Finitistic Dimension~Conjecture.

\begin{dfn}
  Let $(\sfU,\sfV)$ be a complete hereditary cotorsion pair in
  $\catMod{A}$. Let $M$ be an $A$-module.  The
  \emph{$\capUV$-injective dimension of $M$} is defined as
  \begin{equation*}
    \capUV\operatorname{-id}M \deq \inf\left\{ n \in\ZZ \;\middle|\;
      \begin{gathered}
        \text{There is a semi-$\sfU$-$\sfV$ replacement $W$}\\
        \text{of $M$ with $W_{-i}=0$ for $i>n$.}
      \end{gathered}\;\right\} \:.
  \end{equation*}
\end{dfn}

\begin{lem}
  \label{lem:LGor-dimension-modules_b}
  Let $(\sfU,\sfV)$ be a complete hereditary cotorsion pair in
  $\catMod{A}$ and $M$ an $A$-module. There is an inequality
  \begin{equation*}
    \catGL{A}\operatorname{-id}M \dle \capUV\operatorname{-id}{M}
  \end{equation*}
  and equality holds if $\capUV\operatorname{-id}{M}<\infty$.
\end{lem}

\begin{prf*}
  Assume $\capUV\operatorname{-id}M=n$ holds for some integer $n$,
  otherwise the statement is trivial. Let $W$ be a semi-$\sfU$-$\sfV$
  replacement of $M$ with $W_{-i}=0$ for $i>n$. The inequality holds
  as every module in $\capUV$ is left $\sfV$-Gorenstein; see
  \exacite[2.2]{CET-20}. Let $m \le n$ be such that $\H[-i]{M} = 0$
  for $i > m$. If $\Cy[-m]{W}$ is left $\sfV$-Gorenstein, then
  $\Cy[-m]{W}\in \capUV$. Indeed, it has a finite coresolution by
  modules in $\capUV$ which splits by \lemref{sesclosed}(b).
\end{prf*}

\begin{lem}
  \label{lem:LGor-dimension-modules_a}
  Let $(\sfU,\sfV)$ be a complete hereditary cotorsion pair in
  $\catMod{A}$, $M$ an $A$-module in $\sfU$, and $n>0$ an integer. One
  has $\catGL{A}\operatorname{-id} M\leq n$ if and only if there
  exists an exact sequence of $A$-modules
  \begin{equation*}
    0 \lra M \lra G \lra L \lra 0
  \end{equation*}
  with $G\in\catGL{A}$, $L\in \sfU$, and
  $\capUV\operatorname{-id}{L}\leq n-1$.
\end{lem}

\begin{prf*}
  There is, by the completeness of $(\sfU,\sfV)$, a coresolution of
  $M$ consisting of modules in $\capUV$, concentrated in non-positive
  degrees and with cycle modules in $\sfU$; it is a semi-$\sfU$-$\sfV$
  complex by \exaref{bd}. Using this, the ``only if'' statement is
  proved the same way as \thmcite[2.10]{HHl04a}. The ``if'' statement
  holds by a standard application of the Horseshoe Lemma, see
  \lemcite[8.2.1 and Rmk.~8.2.2]{rha}, and \lemref{any-Co}; it also
  follows from the second equality in \rmkref{vanishing}.
\end{prf*}

\begin{prp}
  \label{prp:B-dime-and-GB-dimension}
  Let $(\sfU,\sfV)$ be a complete hereditary cotorsion pair in
  $\catMod{A}$.  For every module $M$ in
  $\sfU\cap {^{\perp}}{\catGL{A}}$ one has
  $\catGL{A}\operatorname{-id}M = \capUV\operatorname{-id}{M}$.
\end{prp}

\begin{prf*}
  By \lemref{LGor-dimension-modules_b}, it suffices to show that
  $\capUV\operatorname{-id}{M} \leq \catGL{A}\operatorname{-id}M$
  holds. Set $\catGL{A}\operatorname{-id}M =n$. If $n=0$, then there
  is an exact sequence of $A$-modules
  \begin{equation*}
    0 \lra H \lra W \lra M \lra 0
  \end{equation*}
  with $H\in\catGL{A}$ and $W \in \capUV$.  Since $M$ is in
  ${^{\perp}{\catGL{A}}}$ the sequence splits, whence $M\in\capUV$. If
  $n\geq1$, then \lemref{LGor-dimension-modules_a} yields an exact
  sequence of $A$-modules, $0 \to M \to G \to L \to 0$, with
  $G\in\catGL{A}$, $L\in{\sfU}$, and
  $\capUV\operatorname{-id}{L}\leq n-1$. There is also an exact
  sequence $0\to G'\to W'\to G\to0$ of $A$-modules with
  $G'\in\catGL{A}$ and $W'\in\capUV$. There is a pullback diagram with
  exact rows and columns:
  \begin{equation*}
    \xymatrix@=1.5pc{&0\ar[d]&0\ar[d]&&\\
      &G'\ar[d]\ar@{=}[r]&G'\ar[d]&&\\
      0\ar[r]&T\ar[r]\ar[d]&W'\ar[r]\ar[d]&L\ar[r]\ar@{=}[d]&0\\
      0\ar[r]&M\ar[r]\ar[d]&G\ar[r]\ar[d]&L\ar[r]&0\\
      &0&\ \, 0\:.}
  \end{equation*}
  As $W'$ is in $\capUV$ and $\capUV\operatorname{-id}{L}\leq n-1$
  holds, exactness of the middle row yields
  $\capUV\operatorname{-id}{T}\leq n$. As the left-hand column splits,
  one has $\capUV\operatorname{-id}{M}\leq n$.
\end{prf*}

Let $\sfP^{<\infty}(A)$ denote the class of finitely generated
$A$-modules of finite projective dimension and $\sfUV$ be the complete
cotorsion pair in $\catMod{A}$ cogenerated by $\sfP^{<\infty}$; that
is, $\sfV = (\sfP^{<\infty}(A))^{\perp}$ and
$\sfU={^{\perp}{\sfV}}$. This cotorsion pair is also hereditary by a
result of Cort\'es Izurdiaga, Estrada, and Guil Asensio
\corcite[2.7]{CEG-12}.

\begin{thm}
  \label{thm:finitistic1}
  Set $\sfV=(\sfP^{<\infty}(A))^{\perp}$ and
  $\sfU={^{\perp}{\sfV}}$. The following conditions are equivalent.
  \begin{eqc}
  \item One has $\findim{A}<\infty$.
  \item Every module in $\sfU$ is $\catGL{A}$-perfect.
  \end{eqc}
\end{thm}

\begin{prf*}
  \proofofimp{i}{ii} If one has $\operatorname{findim}(A)=n$ for an
  integer $n$, then \thmcite[3.4]{CEG-12} yields
  $\capUV\operatorname{-id}{M}\leq n$ for every module $M\in \sfU$; in
  particular $M$ is $\catGL{A}$-perfect by
  \lemref{LGor-dimension-modules_b} and \thmref{perfect}.

  \proofofimp{ii}{i} The coproduct $A^{(A)}$ is a projective
  $A$-module and hence belongs to $\sfU$, so by \thmref{perfect} there
  is an integer $n$ such that $\catGL{A}\operatorname{-id}A^{(A)}=n$
  holds.  As $A^{(A)}$ also belongs to ${^\perp}{\catGL{A}}$, one has
  $\capUV\operatorname{-id}{A^{(A)}}\leq n$ by
  \prpref{B-dime-and-GB-dimension}. Now apply \thmcite[3.4]{CEG-12} to
  conclude that $\findim{A} \leq n$ holds.
\end{prf*}

\begin{cor}
  \label{cor:finitistic}
  Set $\sfV=(\sfP^{<\infty}(A))^{\perp}$ and
  $\sfU={^{\perp}{\sfV}}$. The following conditions are equivalent:
  \begin{eqc}
  \item One has $\operatorname{findim}(A)<\infty$.
  \item The functor
    $\mapdef{\fuFInj}{\stabcatGL{A}}{\catDsgInj{\sfU}}$ yields a
    triangulated equivalence.
  \item One has $\catDdefInj{\sfU}=0$.
  \end{eqc}
\end{cor}

\begin{prf*}
  Combine \thmref{finitistic1} and \corref{A.6}.
\end{prf*}

\begin{rmk}
  Another classic conjecture for Artin algebras is the \emph{Wakamatsu
    tilting conjecture}, which says that any Wakamatsu tilting
  $A$-module of finite projective dimension is a tilting $A$-module;
  see Beligiannis and Reiten \cite[Chap.~IV]{ABlIRt07}. Mantese and
  Reiten \prpcite[4.4]{FMnIRt04} showed that the Wakamatsu tilting
  conjecture is a special case of the Finitistic Dimension
  Conjecture. A Wakamatsu tilting module gives rise to a cotorsion
  pair, see Wang, Li, and Hu \thmcite[1.3]{WLH-19}. If the module has
  finite projective dimension, then vanishing of an associated defect
  category betrays if it is tilting. We omit the details which are
  similar to the arguments~above.
\end{rmk}

\appendix
\stepcounter{section}

\section*{Appendix: A Schanuel's lemma for \semiUV\ replacements}
\label{sec:gordim}

\noindent
Assume throughout this appendix that $\sfA$, as in \secref{gdim}, is
Grothendieck and let $\sfUV$ be a complete hereditary cotorsion pair
in $\sfA$.

\begin{lem}
  \label{lem:2-3}
  Let $0 \to M' \to M \to M'' \to 0$ be an exact sequence in
  $\catC{\sfA}$.
  \begin{prt}
  \item If $M'$ is semi-$\sfV$ then $M$ is semi-$\sfV$ if and only if
    $M''$ is semi-$\sfV$.
  \item If $M''$ is semi-$\sfU$ then $M$ is semi-$\sfU$ if and only if
    $M'$ is semi-$\sfU$.
  \end{prt}
\end{lem}

\begin{prf*}
  We prove part (a); the proof of part (b) is similar. Assume that
  $M'$ is semi-$\sfV$. It is in particular a complex of objects from
  $\sfV$, so $M$ is a complex of objects from $\sfV$ if and only if
  $M''$ is so. Assuming that this is the case, let $U$ be a
  $\sfU$-acyclic complex. As $M'$ is a complex of objects from $\sfV$,
  there is an exact sequence
  \begin{equation*}
    0 \lra \Hom[\sfA]{U}{M'} \lra \Hom[\sfA]{U}{M} \lra
    \Hom[\sfA]{U}{M''} \lra 0\:.
  \end{equation*}
  By assumption the left-hand complex is acyclic, so the middle
  complex is acyclic if and only if the right-hand complex is acyclic.
\end{prf*}

\begin{prp}
  \label{prp:dfn}
  The following assertions hold:
  \begin{prt}
  \item An $\sfA$-complex is $\sfU$-acyclic if and only if it is
    semi-$\sfU$ and acyclic.
  \item An $\sfA$-complex is $\sfV$-acyclic if and only if it is
    semi-$\sfV$ and acyclic.
  \end{prt}
  Moreover, an acyclic \semiUV\ complex is contractible.
\end{prp}

\begin{prf*}
  By \lemcite[3.10]{JGl04} a $\sfU$-acyclic complex is semi-$\sfU$,
  and a $\sfV$-acyclic complex is semi-$\sfV$. The category $\sfA$ has
  enough injectives, see e.g.\ Kashiwara and Schapira
  \thmcite[9.6.2]{cas}, so it follows from \thmcite[3.12]{JGl04} that
  an acyclic semi-$\sfU$ complex is $\sfU$-acyclic. As
  $(\textsf{semi-}\sfU,\sfV\textsf{-ac})$ is complete, see
  \prpcite[3.6]{JGl04} and \prpcite[7.14]{JSt13b}, it follows from
  \lemcite[3.14(1)]{JGl04} that an acyclic semi-$\sfV$ complex is
  $\sfV$-acyclic. Finally, an acyclic \semiUV\ complex is contractible
  as it has cycle objects in $\capUV$.
\end{prf*}

\begin{cor}
  \label{cor:dfn}
  The following assertions hold.
  \begin{prt}
  \item A quasi-isomorphism of semi-$\sfU$ complexes is a
    $\sfU$-quasi-isomorphism.
  \item A quasi-isomorphism of semi-$\sfV$ complexes is a
    $\sfV$-quasi-isomorphism.
  \item A quasi-isomorphism of \semiUV\ complexes is a homotopy
    equivalence.
  \end{prt}
\end{cor}

\begin{prf*}
  The mapping cone of a quasi-isomorphism of semi-$\sfU$ complexes is
  acyclic and by \lemref{2-3} semi-$\sfU$, so it is $\sfU$-acyclic by
  \prpref{dfn}. This proves (a) and the proof of (b) is
  similar. Finally, the mapping cone of a quasi-isomorphism of
  \semiUV\ complexes is an acyclic \semiUV\ complex and hence per
  \prpref[]{dfn} contractible.
\end{prf*}

\begin{lem}
  \label{lem:approx}
  Let $M$ be an $\sfA$-complex and consider the exact sequences from
  \fctref[]{complete},
  \begin{equation*}
    0 \lra V' \lra U \xra{\pi} M \lra 0 \qqand
    0 \lra M \xra{\iota} V \lra U' \lra 0 \:.
  \end{equation*}
  \begin{prt}
  \item If $M$ is a complex of objects in $\sfV$, then $\pi$ is a
    $\sfV$-quasi-isomorphism, and if $M$ is semi-$\sfV$, then $U$ is
    \semiUV.
  \item If $M$ is a complex of objects in $\sfU$, then $\iota$ is a
    $\sfU$-quasi-isomorphism, and if $M$ is semi-$\sfU$, then $V$ is
    \semiUV.
  \end{prt}
\end{lem}

\begin{prf*}
  We prove part (a); the proof of part (b) is similar.  For every
  object $U'$ in $\sfU$ the induced sequence
  \begin{equation*}
    0 \lra \Hom[\sfA]{U'}{V'} \lra \Hom[\sfA]{U'}{U}
    \xra{\Hom[\sfA]{U'}{\pi}} \Hom[\sfA]{U'}{M} \lra 0
  \end{equation*}
  is exact and $\Hom[\sfA]{U'}{V'}$ is acyclic. It follows that
  $\Hom[\sfA]{U'}{\pi}$ is a quasi-isomorphism. If $M$ is a complex of
  objects from $\sfV$, then so is $U$, whence $\Cone{\pi}$ is an
  acyclic complex of objects in $\sfV$. As
  $\Hom[\sfA]{U'}{\Cone{\pi}}$ is acyclic for every $U'\in\sfU$ it
  follows that the cycles of the complex $\Cone{\pi}$ belong to
  $\sfV$. Thus $\pi$ is a $\sfV$-quasi-isomorphism. Finally, if $M$ is
  semi-$\sfV$, then by \prpref{dfn} and \lemref{2-3} so is $U$. That
  is, $U$ is \semiUV.
\end{prf*}

The lemma above provides a \semiUV\ replacement of any $\sfA$-complex;
see \rmkref{replace}. Our next goal is a Schanuel's lemma for \semiUV\
replacements. We move towards it with the next result and prove it in
\prpref{schanuel}.

\begin{prp}
  \label{prp:cothty}
  Let $U$ and $U'$ be semi-$\sfU$ complexes and $V$ and $V'$ be a
  semi-$\sfV$ complexes.
  \begin{prt}
  \item Let $\qisdef{\grb}{U}{U'}$ be a $\sfU$-quasi-isomorphism. For
    every morphism $\mapdef{\gra}{U}{V}$ there is a morphism
    $\mapdef{\grg}{U'}{V}$ with $\grg\grb \sim \gra$, and $\grg$ is
    unique up to homotopy.
  \item Let $\qisdef{\grb}{V'}{V}$ be a $\sfV$-quasi-isomorphism. For
    every morphism $\mapdef{\gra}{U}{V}$ there is a morphism
    $\mapdef{\grg}{U}{V'}$ with $\grb\grg \sim \gra$, and $\grg$ is
    unique up to homotopy.
  \end{prt}
\end{prp}

\begin{prf*}
  We prove part (a); the proof of part (b) is similar. The complex
  $\Cone{\grb}$ is $\sfU$-acyclic, so the morphism
  $\Hom[\sfA]{\grb}{V}$ is a quasi-isomorphism. From this point, the
  proof of \prpcite[1.7]{CELTWY-21} applies verbatim.
\end{prf*}

\begin{lem}
  \label{lem:compare}
  Let $W$ and $W'$ be \semiUV\ complexes. If there is an isomorphism
  $W \qis W'$ in $\catD{\sfA}$, then there is a homotopy equivalence
  $W \to W'$ in $\catC{\sfA}$.
\end{lem}

\begin{prf*}
  By assumption there exists an $\sfA$-complex $X$ and a diagram in
  $\catC{\sfA}$:
  \begin{equation*}
    W \qla X \qra W' \:.
  \end{equation*}
  By \fctref{complete} there is a semi-$\sfU$ complex $U$ and a
  quasi-isomorphism $U \qra X$. The composite $U \lra W$ is by
  \corref{dfn}(a) a $\sfU$-quasi-isomorphism, so up to homotopy the
  composite $U \to W'$ lifts to a quasi-isomorphism $W \lra W'$; see
  \prpref{cothty}. By \corref{dfn}(c) this quasi-isomorphism is a
  homotopy equivalence.
\end{prf*}

\begin{prp}
  \label{prp:schanuel}
  Let $W$ and $W'$ be \semiUV\ complexes isomorphic in
  $\catD{\sfA}$. For every $n\in\ZZ$ the following assertions hold.
  \begin{prt}
  \item There exist objects $U$ and $U'$ in $\sfU\cap \sfV$ such that
    $\Co[n]{W} \oplus U \cong \Co[n]{W'} \oplus U'$.
  \item There exist objects $V$ and $V'$ in $\sfU\cap \sfV$ such that
    $\Cy[n]{W} \oplus V \cong \Cy[n]{W'} \oplus V'$.
  \end{prt}
\end{prp}

\begin{prf*}
  We prove part (a); the proof of part (b) is similar.  By
  \lemref{compare} there is a homotopy equivalence
  $\mapdef{\gra}{W}{W'}$. It follows from \prpref{dfn} that every
  cycle object $\Cy[n]{\Cone{\gra}}$ belongs to $\capUV$. The soft
  truncated morphism $\mapdef{\Tsa{n}{\gra}}{\Tsa{n}{W}}{\Tsa{n}{W'}}$
  is also a homotopy equivalence, so $\Cone{\Tsa{n}{\gra}}$, i.e.\ the
  complex
  \begin{equation*}
    0 \to \Co[n]{W} \to \Co[n]{W'} \oplus W_{n-1}
    \to W'_{n-1} \oplus W_{n-2} \xrightarrow{\dif[n-1]{\Cone{\gra}}}
    W'_{n-2} \oplus W_{n-3} \to \cdots
  \end{equation*}
  is contractible. Hence one has
  $\Co[n]{W} \oplus \Cy[n-1]{\Cone{\gra}} \is \Co[n]{W'} \oplus
  W_{n-1}$.
\end{prf*}

The point of the next result, which is crucial to our proof of
\thmref{perfect}, is that one can exert some control over the
boundedness of \semiUV-replacements of bounded complexes.

\begin{thm}
  \label{thm:bdres}
  The following assertions hold.
  \begin{prt}
  \item For every bounded below semi-$\sfV$ complex $V$ there exists
    an exact sequence
    \begin{equation*}
      0 \lra V' \lra W \xra{\pi} V \lra 0
    \end{equation*}
    such that $V'$ is $\sfV$-acyclic, $W$ is \semiUV\ and bounded
    below, and $\pi$ is a $\sfV$-quasi-isomorphism.
  \item For every bounded above semi-$\sfU$ complex $U$ there exists
    an exact sequence
    \begin{equation*}
      0 \lra U \xra{\iota} W \lra U' \lra 0
    \end{equation*}
    such that $U'$ is $\sfU$-acyclic, $W$ is \semiUV\ and bounded
    above, and $\iota$ is a $\sfU$-quasi-isomorphism.
  \end{prt}
\end{thm}

\begin{prf*}
  We prove part (a); the proof of part (b) is similar. There exists by
  \lemref{approx}(a) an exact sequence
  \begin{equation*}
    0 \lra V' \lra W \xra{\pi} V \lra 0
  \end{equation*}
  such that $V'$ is $\sfV$-acyclic, $W$ is \semiUV\, and $\pi$ is a
  $\sfV$-quasi-isomorphism. Since $V$ is bounded below, the complexes
  $V'$ and $W$ agree in low degrees. It suffices to show that
  $\Cy[n]{W}$ belongs to $\capUV$ for $n \ll 0$, as one then can
  replace $V'$ and $W$ with soft truncations $\Tsb{n}{V'}$ and
  $\Tsb{n}{W}$, cf.~\lemref{2-3} and \exaref{bd}.

  For $n \ll 0$ one has $\Cy[n]{W} = \Cy[n]{V'}$; this is an object in
  $V$ so it remains to show that it is in $\sfU$. Completeness of
  $(\sfU,\sfV)$ provides a resolution $U$ of $\Cy[n]{W}$ that is
  concentrated in non-negative degrees, consists of objects in
  $\capUV$, and has cycle objects in $\sfV$; it is thus a \semiUV\
  replacement per \exaref{bd}.  The truncated complex
  $\Shift[n]{(W_{\le n})}$ is another \semiUV\ replacement of
  $\Cy[n]{W}$, cf.~\lemref{2-3} and \exaref{bd}. Now
  \prpref{schanuel}(b) yields objects $X$ and $X'$ in $\capUV$ such
  that $\Cy[n]{W} \oplus X \is U_0 \oplus X'$ holds. As
  $U_0 \oplus X'$ is in $\sfU$, so is $\Cy[n]{W}$.
\end{prf*}

\section*{Acknowledgment}
\noindent 
We thank the anonymous referee for a very thorough reading and
numerous suggestions that helped improve the paper.



\begin{thebibliography}{10}

\bibitem{MAsMBr69} Maurice Auslander and Mark Bridger, \emph{Stable
    module theory}, Memoirs of the American Mathematical Society,
  No. 94, American Mathematical Society, Providence, R.I.,
  1969. \MR{MR0269685}

\bibitem{rta} Maurice Auslander, Idun Reiten, and Sverre~O. Smal{\o},
  \emph{Representation theory of {A}rtin algebras}, Cambridge Studies
  in Advanced Mathematics, vol.~36, Cambridge University Press,
  Cambridge, 1995. \MR{MR1314422}

\bibitem{HBs60} Hyman Bass, \emph{Finitistic dimension and a
    homological generalization of semi-primary rings},
  Trans. Amer. Math. Soc. \textbf{95} (1960), 466--488.  \MR{MR157984}

\bibitem{ABl00} Apostolos Beligiannis, \emph{The homological theory of
    contravariantly finite subcategories: {A}uslander-{B}uchweitz
    contexts, {G}orenstein categories and (co-)stabilization},
  Comm. Algebra \textbf{28} (2000), no.~10, 4547--4596.
  \MR{MR1780017}

\bibitem{ABlIRt07} Apostolos Beligiannis and Idun Reiten,
  \emph{Homological and homotopical aspects of torsion theories},
  Mem. Amer. Math. Soc. \textbf{188} (2007), no.~883,
  viii+207. \MR{MR2327478}

\bibitem{BJO-15} Petter~Andreas Bergh, Steffen Oppermann, and
  David~A. Jorgensen, \emph{The {G}orenstein defect category},
  Q. J. Math. \textbf{66} (2015), no.~2, 459--471. \MR{MR3356832}

\bibitem{ROB86} Ragnar-Olaf Buchweitz, \emph{Maximal
    {C}ohen--{M}acaulay modules and {T}ate-cohomology over
    {G}orenstein rings}, University of Hannover, 1986, available at
  \mbox{\sffamily http://hdl.handle.net/1807/16682}.

\bibitem{TBh10} Theo B\"{u}hler, \emph{Exact categories},
  Expo. Math. \textbf{28} (2010), no.~1, 1--69. \MR{MR2606234}

\bibitem{TBh11} Theo B\"{u}hler, \emph{On the algebraic foundations of
    bounded cohomology}, Mem. Amer. Math. Soc. \textbf{214} (2011),
  no.~1006, xxii+97. \MR{MR2867320}

\bibitem{CLY-19} Wenjing Chen, Zhongkui Liu, and Xiaoyan Yang, \emph{A
    new method to construct model structures from a cotorsion pair},
  Comm. Algebra \textbf{47} (2019), no.~11, 4420--4431. \MR{MR3991027}

\bibitem{XWC11} Xiao-Wu Chen, \emph{Relative singularity categories
    and {G}orenstein-projective modules}, Math. Nachr. \textbf{284}
  (2011), no.~2-3, 199--212. \MR{MR2790881}

\bibitem{CELTWY-21} Lars~Winther Christensen, Sergio Estrada,
  Li~Liang, Peder Thompson, Dejun Wu, and Gang Yang, \emph{A
    refinement of {G}orenstein flat dimension via the flat-cotorsion
    theory}, J. Algebra \textbf{567} (2021), 346--370.  \MR{MR4159258}

\bibitem{CET} Lars~Winther Christensen, Sergio Estrada, and Peder
  Thompson, \emph{The stable category of {G}orenstein flat sheaves on
    a noetherian scheme}, Proc. Amer.  Math. Soc. \textbf{149} (2021),
  no.~2, 525--538. \MR{MR4198062}

\bibitem{CET-20} Lars~Winther Christensen, Sergio Estrada, and Peder
  Thompson, \emph{Homotopy categories of totally acyclic complexes
    with applications to the flat--cotorsion theory}, Categorical,
  {H}omological and {C}ombinatorial {M}ethods in {A}lgebra,
  Contemp. Math., vol. 751, Amer. Math. Soc., Providence, RI, 2020,
  pp.~99--118. \MR{MR4132086}

\bibitem{LWCPTh19} Lars~Winther Christensen and Peder Thompson,
  \emph{Pure-minimal chain complexes}, Rend. Semin. Mat. Univ. Padova
  \textbf{142} (2019), 41--67.  \MR{MR4032803}

\bibitem{rha} Edgar~E. Enochs and Overtoun M.~G. Jenda, \emph{Relative
    homological algebra}, de Gruyter Expositions in Mathematics,
  vol.~30, Walter de Gruyter \& Co., Berlin, 2000. \MR{MR1753146}

\bibitem{JGl04} James Gillespie, \emph{The flat model structure on
    {${\rm Ch}(R)$}}, Trans.  Amer. Math. Soc. \textbf{356} (2004),
  no.~8, 3369--3390. \MR{MR2052954}

\bibitem{HHl04a} Henrik Holm, \emph{Gorenstein homological
    dimensions}, J. Pure Appl. Algebra \textbf{189} (2004), no.~1-3,
  167--193. \MR{MR2038564}

\bibitem{OImYDn20} Osamu Iyama and Dong Yang, \emph{Quotients of
    triangulated categories and equivalences of {B}uchweitz, {O}rlov,
    and {A}miot-{G}uo-{K}eller}, Amer. J.  Math. \textbf{142} (2020),
  no.~5, 1641--1659. \MR{MR4150654}

\bibitem{SInHKr06} Srikanth Iyengar and Henning Krause,
  \emph{Acyclicity versus total acyclicity for complexes over
    {N}oetherian rings}, Doc. Math. \textbf{11} (2006),
  207--240. \MR{MR2262932}

\bibitem{CEG-12} M.~Cort\'{e}s Izurdiaga, S.~Estrada, and P.~A.~Guil
  Asensio, \emph{A model structure approach to the finitistic
    dimension conjectures}, Math. Nachr.  \textbf{285} (2012), no.~7,
  821--833. \MR{MR2924515}

\bibitem{cas} Masaki Kashiwara and Pierre Schapira, \emph{Categories
    and sheaves}, Grundlehren der Mathematischen Wissenschaften,
  vol. 332, Springer-Verlag, Berlin, 2006. \MR{MR2182076}

\bibitem{BKl90} Bernhard Keller, \emph{Chain complexes and stable
    categories}, Manuscripta Math. \textbf{67} (1990), no.~4,
  379--417. \MR{MR1052551}

\bibitem{BKl96} Bernhard Keller, \emph{Derived categories and their
    uses}, Handbook of algebra, {V}ol. 1, Handb. Algebr., vol.~1,
  Elsevier/North-Holland, Amsterdam, 1996,
  pp.~671--701. \MR{MR1421815}

\bibitem{BKlDVs87} Bernhard Keller and Dieter Vossieck, \emph{Sous les
    cat\'{e}gories d\'{e}riv\'{e}es}, C. R. Acad. Sci. Paris
  S\'{e}r. I Math. \textbf{305} (1987), no.~6,
  225--228. \MR{MR0907948}

\bibitem{HKr10} Henning Krause, \emph{Localization theory for
    triangulated categories}, Triangulated categories, London
  Math. Soc. Lecture Note Ser., vol. 375, Cambridge Univ. Press,
  Cambridge, 2010, pp.~161--235. \MR{MR2681709}

\bibitem{FMnIRt04} Francesca Mantese and Idun Reiten, \emph{Wakamatsu
    tilting modules}, J. Algebra \textbf{278} (2004), no.~2,
  532--552. \MR{MR2071651}

\bibitem{LMaNDn06} Lixin Mao and Nanqing Ding, \emph{The cotorsion
    dimension of modules and rings}, Abelian groups, rings, modules,
  and homological algebra, Lect. Notes Pure Appl. Math., vol. 249,
  Chapman \& Hall/CRC, Boca Raton, FL, 2006,
  pp.~217--233. \MR{MR2229114}

\bibitem{DMf-phd} Daniel Murfet, \emph{The mock homotopy category of
    projectives and {G}rothendieck duality}, Ph.D. thesis, Australian
  National University, 2007, x+145 pp., Available from
  \href{http:www.therisingsea.org}{\sf www.therisingsea.org}.

\bibitem{DMfSSl11} Daniel Murfet and Shokrollah Salarian,
  \emph{Totally acyclic complexes over {N}oetherian schemes},
  Adv. Math. \textbf{226} (2011), no.~2, 1096--1133.  \MR{MR2737778}

\bibitem{ANm90} Amnon Neeman, \emph{The derived category of an exact
    category}, J. Algebra \textbf{135} (1990), no.~2,
  388--394. \MR{MR1080854}

\bibitem{DOO04} D.~O. Orlov, \emph{Triangulated categories of
    singularities and {D}-branes in {L}andau-{G}inzburg models},
  Tr. Mat. Inst. Steklova \textbf{246} (2004),
  no.~Algebr. Geom. Metody, Svyazi i Prilozh.,
  240--262. \MR{MR2101296}

\bibitem{JSt13b} Jan \v{S}\v{t}ov\'{\i}\v{c}ek, \emph{Exact model
    categories, approximation theory, and cohomology of quasi-coherent
    sheaves}, Advances in representation theory of algebras, EMS
  Ser. Congr. Rep., Eur. Math. Soc., Z\"{u}rich, 2013,
  pp.~297--367. \MR{MR3220541}

\bibitem{WLH-19} Jian Wang, Yunxia Li, and Jiangsheng Hu, \emph{When
    the kernel of a complete hereditary cotorsion pair is the additive
    closure of a tilting module}, J.  Algebra \textbf{530} (2019),
  94--113. \MR{MR3938865}

\bibitem{XYnWCh17} Xiaoyan Yang and Wenjing Chen, \emph{Relative
    homological dimensions and {T}ate cohomology of complexes with
    respect to cotorsion pairs}, Comm. Algebra \textbf{45} (2017),
  no.~7, 2875--2888. \MR{MR3594565}
\enlargethispage*{3\baselineskip}
\end{thebibliography}

\def\soft#1{\leavevmode\setbox0=\hbox{h}\dimen7=\ht0\advance \dimen7
  by-1ex\relax\if t#1\relax\rlap{\raise.6\dimen7
    \hbox{\kern.3ex\char'47}}#1\relax\else\if T#1\relax
  \rlap{\raise.5\dimen7\hbox{\kern1.3ex\char'47}}#1\relax \else\if
  d#1\relax\rlap{\raise.5\dimen7\hbox{\kern.9ex
      \char'47}}#1\relax\else\if D#1\relax\rlap{\raise.5\dimen7
    \hbox{\kern1.4ex\char'47}}#1\relax\else\if l#1\relax
  \rlap{\raise.5\dimen7\hbox{\kern.4ex\char'47}}#1\relax \else\if
  L#1\relax\rlap{\raise.5\dimen7\hbox{\kern.7ex
      \char'47}}#1\relax\else\message{accent \string\soft \space #1
    not defined!}#1\relax\fi\fi\fi\fi\fi\fi}
\providecommand{\MR}[1]{\mbox{\href{http://www.ams.org/mathscinet-getitem?mr=#1}{#1}}}
\renewcommand{\MR}[1]{\mbox{\href{http://www.ams.org/mathscinet-getitem?mr=#1}{#1}}}
\providecommand{\arxiv}[2][AC]{\mbox{\href{http://arxiv.org/abs/#2}{\sf
      arXiv:#2 [math.#1]}}} \def\cprime{$'$}
\providecommand{\bysame}{\leavevmode\hbox to3em{\hrulefill}\thinspace}
\providecommand{\MR}{\relax\ifhmode\unskip\space\fi MR }
\providecommand{\MRhref}[2]{%
  \href{http://www.ams.org/mathscinet-getitem?mr=#1}{#2} }
\providecommand{\href}[2]{#2}

\end{document}